\documentclass[reqno]{amsart}
\usepackage[a4paper,centering,vscale=0.75,hscale=0.75]{geometry}
\usepackage[T1]{fontenc}
\usepackage[latin1]{inputenc}
\usepackage[english]{babel}
\usepackage[babel]{csquotes}

\usepackage{cite}
\usepackage{amssymb}
\usepackage{amsmath}
\usepackage{amsthm}
\usepackage{latexsym}
\usepackage{graphicx}
\usepackage{mathrsfs}
\usepackage{mathrsfs}
\usepackage{bbm}
\usepackage{verbatim}

\numberwithin{equation}{section}

\usepackage[usenames,dvipsnames]{color}

\newtheorem{thm}{Theorem}[section]

\newtheorem{defin}[thm]{Definition}
\newtheorem{remark}[thm]{Remark}

\def\enne{\mathbb{N}}

\def\erre{\mathbb{R}}

\def\P{\mathbb{P}}

\def\E{\mathop{{}\mathbb{E}}}

\def\cL{\mathscr{L}}
\def\cF{\mathscr{F}}

\def\eps{\varepsilon}

\def\OO{\mathcal{O}}
\renewcommand{\div}{\operatorname{div}}
\renewcommand{\d}{{\mathrm d}}

\def\beq{\begin{equation}}
\def\eeq{\end{equation}}

\def\to{\rightarrow}
\def\wto{\rightharpoonup}
\def\wstarto{\stackrel{*}{\rightharpoonup}}
\def\embed{\hookrightarrow}
\def\cembed{\stackrel{c}{\hookrightarrow}}

\def\norm #1{\left\|#1\right\|}

\def\sp #1#2{\left<#1,#2\right>}
\newcommand\ip\sp

\begin{document}
\title[Stochastic CH with degenerate mobility and 
logarithmic potential]
{The stochastic Cahn-Hilliard equation\\
with degenerate mobility and logarithmic potential}

\author{Luca Scarpa}
\address{Faculty of Mathematics, University of Vienna
Oskar-Morgenstern-Platz 1, 1090 Vienna, Austria.}
\email{luca.scarpa@univie.ac.at}
\urladdr{http://www.mat.univie.ac.at/$\sim$scarpa}


\begin{abstract}
  We prove existence of martingale solutions for the stochastic 
  Cahn-Hilliard equation with degenerate mobility and
  multiplicative Wiener noise. The potential is allowed to be
  of logarithmic or double-obstacle type.
  By extending to the stochastic framework 
  a regularization procedure 
  introduced by C.~M.~Elliott and H.~Garcke
  in the deterministic setting, 
  we show that a compatibility condition between
  the degeneracy of the mobility
  and the blow-up of the potential
  allows to confine some approximate solutions 
  in the physically relevant domain.
  By using a suitable Lipschitz-continuity property of the noise,
  uniform energy and magnitude estimates are proved.
  The passage to the limit is then carried out 
  by stochastic compactness arguments 
  in a variational framework. Applications to stochastic 
  phase-field modelling are also discussed.
\end{abstract}

\maketitle

%%%%%%%%%%%%%%%%%%%%%%%%%%%%%%%%%%%%%%%%%%%%%%%%

\section{Introduction}
\setcounter{equation}{0}
\label{sec:intro}

The Cahn-Hilliard equation was
firstly proposed in \cite{cahn-hill} in order to describe 
the spinodal decomposition occurring in binary metallic alloys.
Since then it has been increasingly employed 
in several areas such as, among many others,
physics, engineering, and biology.
In the recent years, the Cahn-Hilliard equation 
has become one of the most important 
models involved in phase-field theory. 
In such class of models, the evolution of a certain 
material exhibiting two different features
is described by
introducing a so-called {\em state variable} $\varphi\in[-1,1]$,
representing the difference in volume fractions.
The sets $\{\varphi=1\}$ and $\{\varphi=-1\}$
correspond to the pure regions, while
the interfacial region $\{-1<\varphi<1\}$
where the two components coexist
is supposed to have a positive thickness. 
For this reason, such models are usually referred
to as diffuse interface models, and 
the time evolution of the state variable is often 
described by means of a Cahn-Hilliard-type equation.
The field of applications of diffuse--interface modelling is enormous.
In physics it is used in the context of evolution of separating materials, 
phase--transition phenomena, and dynamics of mixtures of fluids;
in biology phase--field modelling is crucial in the description 
of evolution of interacting cells, tumour growths, and dynamics of interacting 
populations; in engineering it plays a central role 
in modelling of damage and deterioration in continuous media.

Given a smooth bounded domain $\OO$ of $\erre^d$, with $d\geq2$,
and a fixed final time $T>0$,
the deterministic Cahn-Hilliard equation reads
\begin{align}
  \label{eq1_det}
  \partial_t\varphi - \div(m(\varphi)\nabla\mu) = 0 \quad&\text{in } (0,T)\times\OO\,,\\
  \label{eq2_det}
  \mu=-\Delta\varphi + F'(\varphi) \quad&\text{in } (0,T)\times\OO\,,\\
  \label{eq3_det}
  {\bf n}\cdot \nabla\varphi = {\bf n}\cdot m(\varphi)\nabla\mu = 0
  \quad&\text{in } (0,T)\times\partial \OO\,,\\
  \label{eq4_det}
  \varphi(0)=\varphi_0 \quad&\text{in } \OO\,.
\end{align}
The variable $\varphi$ is referred to as state variable,
or order parameter, while $\mu$ is the chemical potential.
Here, the symbol ${\bf n}$ denotes the outward 
unit vector on the boundary $\partial \OO$,
the function $m$ is known as mobility,
while $F:\erre\to[0,+\infty]$ is a double-well potential
with two global minima. Typical examples of $m$ and $F$
are given below. 

The chemical potential $\mu$
is directly related in equation \eqref{eq2_det}
to the subdifferential of the so--called free--energy functional
\[
  \varphi\mapsto\frac12\int_\OO|\nabla\varphi|^2 + \int_\OO F(\varphi)\,.
\]
The double--well potential $F$ may be thought as a singular 
convex function that has been perturbed by a concave quadratic 
function: the effect of the concave perturbation is then 
the creation of two global minima for $F$, each one 
representing the pure phases of the model.
Minimising the $F$--term in
free-energy above describes then the tendency of 
each pure phase to concentrate, whereas 
the gradient term penalises high oscillations of the state variable.
The idea behind the minimisation of the free--energy 
is then a calibration between these two phenomena.

Since only the values of $\varphi$ in the interval $[-1,1]$
are relevant to the physical derivation of the model,
the double--well potential $F$ is only meaningful 
if defined on $[-1,1]$.
The most important example 
is the logarithmic one, defined as
\beq\label{pot_log}
  F_{log}(r):=\frac{\theta}{2}\left((1+r)\ln(1+r) + (1-r)\ln(1-r)\right) + \frac{\theta_0}{2}(1-r^2)\,,
   \quad r\in(-1,1)\,,
\eeq
with $0<\theta<\theta_0$ being given constants, which possesses 
two global minima in the interior of the physically relevant domain $[-1,1]$.
This choice of the nonlinearity is the most coherent with 
the physical derivation of the Cahn--Hilliard model itself, 
in relation to its thermodynamical consistency: for this reason, 
\eqref{pot_log} is usually employed in contexts related to 
separation phenomena in physics.
A second relevant choice for $F$ is the so-called double-obstacle potential
\beq\label{F_ob}
  F_{ob}(r):=
  \begin{cases}
  1-r^2 \quad&\text{if } r\in[-1,1]\,,\\
  +\infty\quad&\text{otherwise}\,.
  \end{cases}
\eeq
Here, by contrast, the global minima corresponds exactly to the pure phases $\pm1$,
and this choice is then often employed in the modelling contexts 
where the pure phases have a privileged role compared to the interface:
it is the cases, for example, of tumor growth dynamics in biology.
In \eqref{F_ob}, the derivative $F_{ob}'$ has to be interpreted 
in the sense of convex analysis as the subdifferential $\partial F_{ob}$,
and the equation takes the form of a differential inclusion.
In some cases, these double-well potentials
are approximated by the polynomial one
\beq\label{F_pol}
  F_{pol}(r):=\frac14(r^2-1)^2\,, \quad r\in\erre\,,
\eeq
which nonetheless does not ensure the relevant constraint $\varphi\in[-1,1]$.
While on the one hand the polynomial potential 
$F_{pol}$ is certainly much easier to handle 
from the mathematical point of view, on the other hand
the logarithmic potential $F_{log}$ is surely the most
relevant in terms of thermodynamical consistency.
Indeed, due to the physical interpretation of diffuse-interface
modelling, only the values of the variable $\varphi$ in $[-1,1]$
are meaningful. For this reason, the possibility 
of dealing with the logarithmic potential $F_{log}$ is crucial.

The classical choice of the mobility $m$ is
a positive constant $m_{con}$, independent of $\varphi$.
Nevertheless, starting from the 
pioneering contribution \cite{cahn-hill2} itself,
several authors proposed 
the choice of a mobility depending 
explicitly on the order parameter.
A thermodynamically relevant choice 
for $m$ has been exhibited in
the works \cite{cahn-tay, cahn-tay2,hill}
and consists of a polynomial mobility $m_{pol}$
defined on the physically relevant domain $[-1,1]$
with degeneracy at the extremal points:
\beq\label{m_pol}
  m_{pol}(r):=1-r^2\,, \quad r\in[-1,1]\,.
\eeq
A more general version of $m_{pol}$ is given by
\beq\label{m_alpha}
  m_{\alpha}(r):=(1-r^2)^\alpha\,, \quad r\in[-1,1]\,, \quad \alpha\geq1\,.
\eeq

Several variants of the Cahn-Hilliard equation have been
studied in the last decades. 
Novick-Cohen proposed in \cite{novick-cohen} the 
viscous regularization (see also \cite{ell-st, ell-zhen}), 
accounting also for viscous dynamics 
occurring in phase-transition evolution.
Gurtin generalized the viscous correction 
in \cite{gurtin} possibly including 
nonlinear viscosity contributions in the 
equation.
More recently, physicists have introduced 
the so-called dynamic boundary conditions
in order to account also for possible interaction
with the walls in a confined system (see
for example \cite{fish-spinod, kenz-spinod, gal-DBC2}).

The mathematical literature on the deterministic 
Cahn-Hilliard equation is extremely developed:
we refer to \cite{mir-CH} and the references therein for 
a unifying treatment on the available literature.
In particular, in the case of constant mobility
existence, uniqueness, and regularity 
have been studied in
\cite{col-fuk-eqCH, col-gil-spr, cher-gat-mir, 
cher-mir-zel, colli-fuk-CHmass,
gil-mir-sch} both with irregular potentials and
dynamic boundary conditions,
and in \cite{bcst1, bcst2, mir-sch, scar-VCHDBC}
with nonlinear viscosity terms.
Significant attention has been devoted also 
to the asymptotic behaviour of solutions
\cite{col-fuk-diffusion, col-scar, gil-mir-sch-longtime}
and optimal control problems
\cite{col-far-hass-gil-spr, col-gil-spr-contr, col-gil-spr-contr2, hinter-weg}.
A mathematical analysis of the framework
of nonconstant and possibly degenerate mobility
has been investigated in \cite{cahn-ell-nov,ell-gar}.
In this direction we also refer to the contribution
\cite{chang, dalp-giac-nov-ACHDM, jing}
regarding existence of solutions.
Let us point out the work
\cite{grass-mir-ros-sch-CHDM} dealing with
the analysis of a Cahn-Hilliard equation
with mobility depending on the chemical potential,
\cite{sch} in relation to global attractors,
and \cite{lis-matt-sav-CHDM} for an approach
based on gradient flows in Wasserstein spaces.
A diffuse interface model with degenerate mobility 
has been studied also in \cite{frig-grass-roc2}.
Numerical simulations have been analyzed in
\cite{barr-blow-gar-CHDM, kim}.

Despite the fact that the deterministic model has been proven 
to be extremely effective in the description of phase--separation, 
there are certainly important downsides.
One of the main drawbacks of the deterministic framework is the
impossibility of describing the unpredictable disruptions 
occurring in the evolution at the microscopic scale.
These may be due to several phenomena of different nature,
such as uncertain movements at a microscopic level
caused by configurational, electronic, or magnetic effects,
which cannot be captured 
by the classical deterministic Cahn--Hilliard system.
The most natural way to capture the randomness 
component which may affect phase-field evolutions
is to introduce a Wiener-type noise in the 
Cahn-Hilliard equation itself.
This was first proposed in \cite{cook} 
employing Wiener noises in the well--known 
Cahn--Hilliard--Cook stochastic version of the model,
which has been then validated as the 
only genuine description of phase--separation
in the contributions \cite{bin-SCH, pego-SCH}.
Currently, this is widely studied both in the physics and 
applied mathematics literatures, for which we refer to 
\cite{rer-SCH, mhb-SCH}
and the references therein.
The stochastic version Cahn-Hilliard equation reads
\begin{align}
  \label{eq1}
  \d \varphi - \div\left(m(\varphi)\nabla\mu\right)\,\d t = G(\varphi)\,\d W \qquad&\text{in } (0,T)\times\OO\,,\\
  \label{eq2}
  \mu = -\Delta\varphi + F'(\varphi) \qquad&\text{in } (0,T)\times\OO\,,\\
  \label{eq3}
  {\bf n}\cdot\nabla\varphi = {\bf n}\cdot m(\varphi)\nabla\mu= 0 
  \qquad&\text{in } (0,T)\times\partial\OO\,,\\
  \label{eq4}
  \varphi(0)=\varphi_0 \qquad&\text{in } \OO\,.
\end{align}
Here, $W$ is a cylindrical Wiener process defined on 
certain stochastic basis, while $G$ is a suitable 
stochastically integrable operator with respect to $W$.
Precise assumptions on the data are given in Section~\ref{sec:main} below.

From the mathematical point of view, 
the stochastic Cahn-Hilliard equation has been
studied mainly in the case of polynomial potential
and only with constant mobility.
One the first contributions in this direction
is \cite{daprato-deb}, in which the authors show
existence of solutions via a semigroup approach
in the case of polynomial potentials.
More recently, well-posedness has been 
investigated also in \cite{corn, elez-mike}
again in the polynomial setting.
A more general framework allowing 
for rapidly growing potentials (e.g.~more than exponentially) 
has been analyzed in \cite{scar-SCH, scar-SVCH}
from a variational approach.
The genuine case of logarithmic potentials
has only been covered in the works \cite{deb-zamb, deb-goud, goud}
by means of so-called reflection measures.

The mathematical literature on stochastic phase-field modelling
has also been increasingly developed.
Let us point out in this direction the works
\cite{ant-kar-mill} dealing with unbounded noise,
\cite{feir-petc,feir-petc2} for a study of a diffuse interface model 
with thermal fluctuations,
and \cite{bauz-bon-leb, orr-scar}
dealing with the stochastic Allen-Cahn equation. 
Beside well-posedness,
optimal control problems have also
been studied in \cite{scar-OCSHC}
in the case of the
stochastic Cahn-Hilliard equation, and in
\cite{orr-roc-scar} in the context of a 
stochastic phase-field model for tumour growth.
Besides, let us point out the mathematical literature 
on stochastic two--phase flows 
has also been expanded in the last years, in the context of 
coupled stochastic systems of Cahn--Hilliard--Navier--Stokes 
and Allen--Cahn--Navier--Stokes type.
In this direction, we refer the reader to the contributions 
\cite{DTM-stochCHNS4, DTM-stochCHNS6, TM-stoch2phase, 
TM-stochCHNS} for existence of solutions,
\cite{barbu-stoch2phase, DTM-stochCHNS5} about asymptotic long--time behaviour,
\cite{DTM-stochCHNS3} on large deviation limits, 
and \cite{DTM-stochCHNS1, DTM-stochCHNS2} dealing
with a nonlocal phase--field equation in the system instead.

One of the most critical problems of the stochastic model
is that the presence of the random forcing term does not guarantee 
that the state variable $\varphi$ remains in the physically relevant domain $[-1,1]$,
even if the double--well potential is singular as in \eqref{pot_log} or \eqref{F_ob}.
This is due to the presence of the additional second--order term in the energy balance,
which may cause blow--up of the energy in finite time.
The consequences of this fact are sever, both on the mathematical side and
especially from the modelling perspective. Indeed, 
from the point of view of thermodynamically consistency of the equation, 
$\varphi$ represents a local concentration and is only meaningful 
if belonging to the physical interval $[-1,1]$: the impossibility of 
proving that the solution of the stochastic equation satisfies this constraint 
inevitably represents a modelling downside.
Besides, as we have pointed out above, the available literature on the stochastic 
Cahn--Hilliard equation only handles the case of constant mobility,
which unfortunately is not the most suited for describing phase--separation.

The current state--of--the--art of the stochastic Cahn--Hilliard equation 
inevitably calls then for a deeper investigation in the direction of 
degeneracy of the mobility, including the physically relevant logarithmic potential,
and showing that the physically meaningful constraint $\varphi\in[-1,1]$ is achieved.
This paper is the first contribution in the mathematical literature 
that addresses these three points, and represents then 
an important step especially in terms of applicative validation of the
stochastic Cahn--Hilliard--Cook model.
Our results have important consequences 
to all fields, in particular physics and engineering, where
phase--separation is usually studied under random forcing,
as we provide the first mathematical validation of the stochastic model 
in its more relevant form, i.e.~with degenerate mobility $m_{pol}$
and logarithmic potential $F_{log}$.

In this paper we are interested in studying the stochastic 
Cahn-Hilliard equation \eqref{eq1}--\eqref{eq4}
from a variational approach, including the cases
of degenerate mobility $m_{pol}$, the logarithmic 
double-well potential $F_{log}$, and the double-obstacle potential $F_{ob}$. 
As we have pointed out, these choices are indeed 
the most relevant in terms of the thermodynamical 
coherency of the model. 
From the mathematical perspective, our approach 
differs from \cite{deb-zamb, deb-goud, goud} as we do not
rely on reflection measures.
Our techniques extend to stochastic framework 
the ideas of C.~M.~Elliott and H.~Garcke in \cite{ell-gar},
and consist of a compatibility condition between 
the degeneracy of the mobility, the
coefficient $G$, and the possible blow-up of
the potential at the extremal points $\pm1$.

Let us briefly explain the main difficulties arising in the case of 
degenerate mobility and logarithmic potential, and how 
we overcome these in the present work.

The first main issue appearing in the stochastic framework 
is the presence of a proliferation term in equation \eqref{eq1}.
Indeed, in the deterministic case
the integration in space of equation \eqref{eq1_det}
yields, together with the boundary conditions \eqref{eq3_det},
the conservation of mass during the evolution. 
This is in turn crucial when dealing with irregular potentials
such as $F_{log}$ or $F_{ob}$, as it allows 
to control the spatial mean of the chemical potential.
However, in the stochastic scenario
the presence of the noise term in equation \eqref{eq1}
determines a proliferation of the total mass of the system.
Whereas this drawback can be 
overcome in the easier case of regular potentials as $F_{pol}$,
this results in the impossibility of obtaining 
satisfactory estimates on the chemical potential $\mu$
in the case of logarithmic and double-obstacle potentials.
The main reason is the following.
On the one hand the derivatives of $F_{pol}$
can be controlled by $F_{pol}$ itself, i.e.
$|F''_{pol}|,\, |F'_{pol}| \leq c(1+ F)$
for a certain $c>0$, so that the usual energy estimates on $F_{pol}$
allow to bound also $F'_{pol}$, hence $\mu$ as a byproduct.
On the other hand, however, the derivatives of the logarithmic potential
$F_{log}$ blow up at $\pm1$ much more rapidly than $F_{log}$ itself,
so that the classical energy estimates on $F_{log}$ are not enough 
to deduce a control on $F_{log}'$.
This problem is even more evident
in the stochastic setting due to
the presence in the energy estimates
of the second order It\^o correction, which depends of the 
second derivative $F''$. 
Again, while this term can be handled in 
the case of polynomial potentials as $F_{pol}$,
in the case of logarithmic potential $F_{log}$
the situation is much more critical, due to 
blow-ups at $\pm1$ pointed out above.

The second main issue is
the degeneracy of the mobility. Indeed, while in the case 
of constant mobility $m_{con}$, or more generally if $m$ is bounded
from below by a positive constant, one usually deduces 
estimates on $\nabla\mu$ pretty directly, if $m$
degenerates at $\pm1$, as in the physically relevant case $m_{pol}$,
then there is no hope to obtain a control on the gradient of the chemical potential.

These two main problems suggest
that in the case of degenerate mobility 
the role of the chemical potential $\mu$ must be 
passed by, and a different interpretation of the equation is needed.
To this end, if one formally substitutes equation \eqref{eq2} into \eqref{eq1},
it is possible to obtain a variational formulation on the problem 
only involving the variable $\varphi$. In particular, the nonlinear
term resulting from such substitution (see Definition~\ref{def:sol-D} below)
is in the form $m(\varphi)F''(\varphi)$. Hence, supposing that 
the degeneracy of the mobility compensates the blow-up of $F''$
at $\pm1$, we can obtain a coherent formulation of the problem
not involving $\mu$ anymore.
Such compatibility condition between $m$ and $F''$
was employed in the mentioned work \cite{ell-gar}, and
is very natural as it is satisfied by the physically relevant choices $m_{pol}$
and $F_{log}$, as we will show in Remark~\ref{ex:log} below.

The idea to overcome the presence of It\^o correction 
terms in the energy estimates depending on $F''$ is of similar nature.
If $G$ is Lipschitz-continuous and vanishes at the extremal points $\pm1$,
its degeneracy can compensate the blow-up of $F''$ and 
the energy estimate can be closed. Again, such Lipschitz-continuity 
assumption on $G$ is very natural in applications, 
and has been widely employed
in stochastic phase-field modelling.
For example, in \cite{bauz-bon-leb} the authors 
consider a diffuse--interface model based on the 
stochastic Allen--Cahn equation in the context of 
evolution of damage in certain continuum media. 
The state variable represents here the damage parameter,
which may be though as the local ratio of active cohesive bonds
at the microscopic level. Clearly, 
only positive values of the state variable are physically meaningful,
and in order to ensure this the authors use a 
obstacle--type potential. As far as the noise is concerned, 
the idea to handle the singularity is of similar nature:
the random forcing is supposed to be multiplicative and 
to ``switch off'' whenever the $\varphi$ touches the potential barriers.  
On the same line, the superposition--type stochastic forcing 
that we propose in this paper has also been widely 
employed in applications. In \cite{feir-petc} for instance,
this has been actively explored in the context of 
a model for a binary mixture of incompressible fluids.

In this work we prove two main results. The first one deals
with existence of martingale solutions to the problem \eqref{eq1}--\eqref{eq4}
in the case of positive mobility and regular potential.
Although being a preparatory work for us, this first
result is also interesting on its own, as it covers the case of 
non constant mobility and allows up to first-order 
exponential growth for the potential. The proofs rely on a double 
approximation involving a Faedo-Galerkin discretization
in space and an Yosida regularization on the nonlinearity.
The second main result that we prove is the 
main contribution of the paper, and states 
existence of martingale solutions in the 
case of degenerate mobility and irregular potentials,
possibly including $F_{log}$ and $F_{ob}$. We
employ a suitable regularisation on the potential,
the mobility, and the diffusion coefficient so that we can 
solve the approximated problem thanks to the first result.
We show then uniform estimates on the solutions
by using energy and magnitude estimates based on the 
compatibility between $m$, $F$, and $G$. Finally we pass to 
the limit by a stochastic compactness argument.

As far as uniqueness of solutions is concerned, 
the problem is still open even in the deterministic setting.
A uniqueness result for the system with degenerate mobility 
has been obtained in \cite{frig-lam-roc} in the framework 
of a nonlocal diffusion related to tumour growth dynamics:
here, the authors exploit the regularizing properties of
the nonlocal nature of the equation in order to
show continuous dependence on the initial data.
Nevertheless, for Cahn-Hilliard evolutions of 
local type with degenerate mobility, regularity 
of solutions is much more difficult to achieve, and 
uniqueness remains unknown. In the stochastic 
setting, this also prevents from proving 
existence of probabilistically strong solutions.

Let us now summarize the main contents of the work.
In Section~\ref{sec:main} we introduce the main setting 
and we state the main results. Section~\ref{sec:ND}
contains the proof of existence of martingale solutions
in the case of positive mobility and regular potential.
Section~\ref{sec:D} is focused on the proof of existence 
of martingale solutions in the setting of degenerate mobility 
and irregular potential.

%%%%%%%%%%%%%%%%%%%%%%%%%%%%%%%%%%%%%%%%%%%%%%%%

\section{Main results}
\label{sec:main}

We introduce here the notation and setting of the paper,
and state the main results. The first main result focuses
on existence of solutions in case of nondegenerate mobility 
and regular potential, while the second deals with the 
case of degenerate mobility and logarithmic potential.

\subsection{Notation and setting}
\label{ssec:not}

For any real Banach space $E$, its dual will be denoted by $E^*$.
The norm in $E$ and the duality pairing between $E^*$ and $E$ will be denoted by 
$\norm{\cdot}_E$ and $\ip{\cdot}{\cdot}_E$, respectively.
If $(A,\mathscr{A},\nu)$ is a finite measure space, 
we use the classical notation $L^p(A; E)$ for the space
of $p$-Bochner integrable functions, for any $p\in[1,+\infty]$.
We shall also use the classical symbol $L^0(A; E)$ for the space
of $\mathscr{A}$-measurable functions with values in $E$.
If $E_1$ and $E_2$ are separable Hilbert spaces, 
we use the notation $\cL^2(E_1,E_2)$ for the space of 
Hilbert-Schmidt operators from $E_1$ to $E_2$.

Throughout the paper, $(\Omega,\cF, (\cF_t)_{t\in[0,T]}, \P)$ is a filtered probability space
satisfying the usual conditions, with $T>0$ being a fixed final time,
and $W$ is a cylindrical Wiener process on a separable Hilbert space $U$.
We fix once and for all a complete orthonormal system $(u_k)_k$ of $U$.
For every separable Hilbert space $E$ and $\ell\in\mathopen[2,+\infty\mathclose)$, we set
\[
  L^\ell_w(\Omega; L^\infty(0,T; E^*)):=\left\{v:\Omega\to L^\infty(0,T; E^*) \text{ weak*-measurable}:\;
  \E\norm{v}_{L^\infty(0,T; E^*)}^\ell<+\infty
  \right\}\,,
\]
and recall that by \cite[Thm.~8.20.3]{edwards} we have
\[
L^\ell_w(\Omega; L^\infty(0,T; E^*))=
\left(L^{\frac{\ell}{\ell-1}}(\Omega; L^1(0,T; E))\right)^*\,.
\]
Moreover, we will use the symbols $C^0([0,T]; E)$ and $C^0_w([0,T]; E)$
for the spaces of strongly and weakly continuous functions from $[0,T]$ to $E$,
respectively.

  It is useful to recall here some general facts about the cylindrical Wiener process $W$
  and stochastic integration that will be used later on in the paper:
  we follow the approach of \cite[\S 2.5.1--2]{LiuRo}.
  Since $W$ is cylindrical in $U$, we have the formal representation 
  \beq
  \label{W_cil}
  W=\sum_{k=0}^\infty\beta_ku_k\,,
  \eeq
  with $(\beta_k)_k$ being a family of independent real Brownian motions.
  Nonetheless, it is important to note that the infinite sum does {\em not} converge 
  in general in $U$, hence $W$ is not rigorously defined as a $U$--valued continuous process.
  In order to properly define $W$ and the corresponding stochastic integral,
  it is useful to show that it is always possible to consider $W$ as a 
  $Q_1$--Wiener process on a larger space $U_1$, with 
  $Q_1$ being of trace--class on $U_1$.
  To this end, note that there always exists a larger separable Hilbert space 
  $U_1$ and a Hilbert--Schmidt operator $\iota\in\cL^2(U,U_1)$. For example, 
  on $U$ one can define the norm
  \[
  \norm{v}_{U_1}:=\left(\sum_{k=0}^\infty\frac{1}{k^2}|(v,u_k)_U|^2\right)^{1/2}\,, \qquad v\in U\,.
  \]
  It is easy to check that $\norm{\cdot}_{U_1}$ is a well--defined norm on $U$, 
  weaker than the usual one $\norm{\cdot}_U$. It makes sense then to 
  define $U_1$ as the abstract closure of $U$ with respect to 
  the norm $\norm{\cdot}_{U_1}$:
  namely, $U_1$ is the abstract space of infinite linear combinations $v$ of $(u_k)_k$
  for which $\norm{v}_{U_1}$ is finite, i.e.
  \[
  U_1:=\left\{v=\sum_{k=0}^\infty\alpha_ku_k:
  \quad
  \norm{v}_{U_1}^2:=\sum_{k=0}^\infty\frac{\alpha_k^2}{k^2}<+\infty
  \right\}\,.
  \]
  With this choice, $(U_1, \norm{\cdot}_{U_1})$ is actually a separable Hilbert space 
  and the inclusion $\iota:U\embed U_1$ is Hilbert--Schmidt: indeed,
  we have 
  \[
  \sum_{k=0}^\infty\norm{\iota (u_k)}_{U_1}^2=\sum_{k=0}^\infty\norm{u_k}_{U_1}^2
  =\sum_{k=0}^\infty\sum_{j=0}^\infty\frac{1}{j^2}|(u_k,u_j)_U|^2
  =\sum_{k=0}^\infty\frac1{k^2}\norm{u_k}_U^2=\sum_{k=0}^\infty\frac1{k^2}<+\infty\,.
  \]
 Consequently, by the properties of Hilbert--Schmidt operators (see again \cite{LiuRo}),
 we have that the infinite sum \eqref{W_cil} actually converges in $U_1$: this means 
 that we can look at $W$ 
 as a rigorously defined stochastic process on $U_1$. Moreover, it actually holds that 
 $W$ is a well--defined $Q_1$--Wiener process on $U_1$, 
 with $Q_1:=\iota\circ\iota^*$ being of trace class on $U_1$,
 and such that $Q_1^{1/2}(U_1)=\iota(U)$.
 In the sequel we will say that
 $W$ is a cylindrical Wiener process on $U$ if and only if
 it is a $Q_1$--Wiener process on $U_1$.
 Furthermore, 
 stochastic integration with respect to the cylindrical process $W$
 is defined in terms of the usual stochastic integration 
 with respect to the $Q_1$--Wiener process.
 In this regard, 
 for every Hilbert space $K$ it holds that $B\in\cL^2(U,K)$ if and only if
  $B\circ \iota^{-1}\in\cL^2(Q_1^{1/2}(U_1), H)$.
 Consequently, 
 for every progressively measurable stochastic integrand 
 $B\in L^2(\Omega; L^2(0,T; \cL^2(U,K)))$, it holds that 
 $B\circ \iota^{-1}\in L^2(\Omega; L^2(0,T; \cL^2(Q_1^{1/2}(U_1),K)))$ and
 \[
 \int_0^\cdot B(s)\,\d W(s) := \int_0^\cdot B\circ \iota^{-1}(s)\,\d W(s)\,,
 \] 
 where the right--hand side is the usual integral with respect to the $Q_1$--Wiener 
 process, the left--hand side is the stochastic integral with respect to the cylindrical 
 Wiener process, and the equality is intended in the sense of indistinguishable 
 continuous $K$--valued processes. The definition of stochastic integral with respect
 to the cylindrical Wiener process $W$ can be shown to be independent of
 the specific Hilbert space $U_1$ and the Hilbert--Schmidt embedding $\iota$.

Let $\OO\subset\erre^d$ ($d\geq2$) be a smooth bounded domain.
We shall use the notation $Q:=(0,T)\times\OO$ and $Q_t:=(0,t)\times\OO$
for every $t\in(0,T)$.
Denoting by ${\bf n}$ the outward normal unit vector on $\OO$,
we define the functional spaces 
\[
  H:=L^2(\OO)\,,  \qquad
  V_1:=H^1(\OO)\,, \qquad
  V_2:=\{v\in H^2(\OO):\;{\bf n}\cdot\nabla v = 0\quad\text{a.e.~on } \partial\OO\}\,,
\]
endowed with their natural norms $\norm{\cdot}_H$, 
$\norm{\cdot}_{V_1}$, and $\norm{\cdot}_{V_2}$,
respectively. For every $v\in V_1^*$, we set $v_\OO:=\frac1{|\OO|}\ip{v}{1}$
for the spatial mean of $v$. We also define
\[  
\mathcal B_R:=\{v\in L^\infty(\OO): \norm{v}_{L^\infty(\OO)}\leq R\}\,, \quad R>0\,.
\]

Moreover, we will use the symbol $c$ to denote any arbitrary positive
constant depending only on the data of the problem, whose value
may be updated throughout the proofs. When we want to specify 
the dependence of $c$ on specific quantities, we will indicate them
through a subscript.

\subsection{Nondegenerate mobility and regular potential}
In case of nondegenerate mobility and regular potential,
we assume the following.
\begin{itemize}
  \item[{\bf ND1}] $F\in C^2(\erre)$, $F\geq0$, $F'(0)=0$, and there exists a constant $C_F>0$ such that 
  \begin{align*}
  |F'(r)| \leq C_F\left(1 + F(r)\right) \qquad&\forall\,r\in\erre\,, \\
  |F''(r)| \leq C_F\left(1 + F(r)\right) \qquad&\forall\,r\in\erre\,, \\
  F''(r)\geq -C_F \qquad&\forall\,r\in\erre\,.
  \end{align*}
  \item[{\bf ND2}] $m\in C^0(\erre)$ and there exist two constants $m_*,m^*>0$ such that 
  \[
  m_*\leq m(r)\leq m^* \quad\forall\,r\in\erre\,.
  \]
  \item[{\bf ND3}] $G: H\to \cL^2(U,H)$ is measurable and there exists 
  $(g_k)_{k\in\enne}\subset W^{1,\infty}(\erre)$ such that
  \begin{align*}
  &G(v)u_k=g_k(v) \qquad\forall\,v\in H,\quad\forall\,k\in\enne\,,\\
  &C_G:=\sum_{k=0}^\infty\norm{g_k}^2_{W^{1,\infty}(\erre)}<+\infty\,.
  \end{align*}
  \item[{\bf ND4}] The initial datum in nonrandom, and satisfies
  $\varphi_0\in V_1$ and $F(\varphi_0)\in L^1(\OO)$.
\end{itemize}
Let us point out that the assumption {\bf ND1} allows for the classical choice
of the polynomial double-well potential $F_{pol}$ defined in \eqref{F_pol},
but also allows to consider polynomials of any orders
and even first-order exponentials. Moreover, 
assumption {\bf ND2} allows of course for the
constant mobility scenario, but also
includes the case of positive nonconstant 
mobilities. Condition {\bf ND3} on the noise
is widely employed in literature
(see for example
\cite{breit-feir-hof,feir-petc, feir-petc2}), and ensures that
in particular that 
$G:H\to\cL^2(U,H)$ is Lipschitz-continuous and linearly bounded, 
and that the restriction $G_{|V_1}:V_1\to \cL^2(U,V_1)$ is linearly bounded.
The initial datum is assumed to be nonrandom in {\bf ND4}:
this is meaningful in relation to the physical interpretation of the model.
A random initial datum could also be considered, but this would 
make the mathematical treatment too heavy in our opinion,
as different estimates are based on different moments in $\Omega$.
Since this is not the main focus of the paper, we preferred to 
assume $\varphi_0$ nonrandom in order to make the treatment clearer
and avoid technicalities: for an exact analysis on the moments of the 
initial datum we refer the reader to \cite{scar-SVCH}.

We precise now the definition of martingale solution
in the case of nondegenerate mobility and regular potential.
\begin{defin}\label{def:sol-ND}
  Assume conditions {\bf ND1}--{\bf ND4}.
  A martingale solution to the problem \eqref{eq1}--\eqref{eq4} is a septuple
  $(\hat\Omega, \hat\cF, (\hat\cF_t)_{t\in[0,T]}, \hat\P, \hat W, \hat\varphi, \hat\mu)$
  such that:
  \begin{itemize}
  \item $(\hat\Omega, \hat\cF, (\hat\cF_t)_{t\in[0,T]}, \hat\P)$ is a filtered probability 
  space satisfying the usual conditions;
  \item $\hat W$ is a $U$-valued cylindrical Wiener process on 
  $(\hat\Omega, \hat\cF, (\hat\cF_t)_{t\in[0,T]}, \hat\P)$;
  \item $\hat\varphi \in L^0(\hat\Omega; C^0([0,T]; H)\cap L^2(0,T; V_2))$
  is progressively measurable;
  \item $\hat\mu=-\Delta\hat\varphi + F'(\hat\varphi) \in L^0(\hat\Omega; L^2(0,T; V_1))$;
  \item for every $v\in V_1$, it holds that 
  \beq\label{var-for-ND}\begin{split}
  \int_\OO\hat\varphi(t,x)v(x)\,\d x &+ 
  \int_{Q_t}m(\hat\varphi(s,x))\nabla\hat\mu(s,x)\cdot\nabla v(x)\,\d x\,\d s \\
  &=\int_\OO\varphi_0(x)v(x)\,\d x + 
  \int_\OO\left(\int_0^tG(\hat\varphi(s))\,\d \hat W(s)\right)\!\!(x)v(x)\,\d x
  \end{split}\eeq
  for every $t\in[0,T]$, $\hat\P$-almost surely.
  \end{itemize}
\end{defin}

The first main result of the paper deals with existence of martingale solutions
in case of positive mobility and regular potential.

\begin{thm}
  \label{thm1}
  Assume conditions {\bf ND1}--{\bf ND4}. Then, 
  there exists a martingale solution
  \[
  \left(\hat\Omega, \hat\cF, (\hat\cF_t)_{t\in[0,T]}, \hat\P, \hat W, \hat\varphi, \hat\mu\right)
  \]
  to the problem \eqref{eq1}--\eqref{eq4} in the sense of Definition~\ref{def:sol-ND}
  such that, for every $\ell\in\mathopen[2,+\infty\mathclose)$,
  \begin{align*}
  &\hat\varphi \in L^\ell(\hat\Omega; C^0([0,T]; H))\cap 
  L^\ell_w(\hat\Omega;L^\infty(0,T; V_1))\cap L^{\ell/2}(\hat\Omega; L^2(0,T; V_2))\,,\\
  &\hat\mu \in L^{\ell/2}(\hat\Omega; L^2(0,T; V_1))\,, \qquad
  \nabla\hat\mu\in L^\ell(\hat\Omega; L^2(0,T; H^d))\,,\\
  &F'(\hat\varphi) \in L^{\ell/2}(\hat\Omega; L^2(0,T; H))\,,
  \end{align*}
  and the following energy inequality holds, for every $t\in[0,T]$:
  \beq\label{en_ineq}
  \begin{split}
    \frac12\sup_{r\in[0,t]}\hat\E\norm{\nabla\hat\varphi(r)}_H^2
    &+\sup_{r\in[0,t]}\hat\E\norm{F(\hat\varphi(r))}_{L^1(\OO)}
    +\hat\E\int_{Q_t}m(\hat\varphi(s,x))|\nabla\hat\mu(s,x)|^2\,\d x\,\d s\\
    &\leq\frac12\norm{\nabla\varphi_0}_H^2 + \norm{F(\varphi_0)}_{L^1(\OO)}
    +\frac{C_G}{2}\hat\E\int_0^t\norm{\nabla\hat\varphi(s)}_H^2\,\d s\\
    &+\frac12\hat\E\int_0^t\sum_{k=0}^\infty\int_\OO
    F''(\hat\varphi(s,x))|g_k(\hat\varphi(s,x))|^2\,\d x\,\d s\,.
    \end{split}
  \eeq
  If also 
  \beq\label{extra_growth}
  |F''(r)|\leq C_F(1+|r|^q) \quad\forall\,r\in\erre\,, \qquad\text{where}\quad
  \begin{cases}
  q\in\mathopen[2,+\infty\mathclose) \quad&\text{if } d=2\,,\\
  q:=\frac2{d-2} \quad&\text{if } d\geq 3\,,
  \end{cases}
  \eeq
  then it holds that
  \[
  \hat\varphi\in L^{\ell/4}(\hat\Omega; L^2(0,T; H^3(\OO)))\,, \qquad
  F'(\hat\varphi)\in L^{\ell/4}(\hat\Omega; L^2(0,T; V_1))\,.
  \]
\end{thm}

\subsection{Degenerate mobility and irregular potential}
We deal now with a degenerate mobility $m:[-1,1]\to\erre$ which
vanishes at $\pm 1$ and an irregular potential $F:(-1,1)\to\erre$
possibly of logarithmic or double-obstacle type. 
In this case we assume the following.

\begin{itemize}
  \item[{\bf D1}] $F:(-1,1)\to\mathopen[0,+\infty\mathclose)$ can be decomposed as
  $F=F_1+F_2$, where $F_1\in C^2(-1,1)$ is convex and $F_2\in C^2([-1,1])$.
  \item[{\bf D2}] $m\in W^{1,\infty}(-1,1)$ is such that
  \[
  m(r)\geq0 \quad\forall\,r\in[-1,1]\,, \qquad
  m(r)=0 \quad\text{iff } r=\pm1\,, \qquad
  mF''\in C^0([-1,1])\,.
  \]
  In particular, it is well defined the function
  \[
  M:(-1,1)\to[0,+\infty)\,, \qquad
  M(0)=M'(0)=0\,, \qquad
  M''(r)=\frac1{m(r)}\,, \quad r\in(-1,1)\,.
  \]
  \item[{\bf D3}] $G: \mathcal B_1\to \cL^2(U,H)$ 
  is measurable and there exists
  $(g_k)_{k\in\enne}\subset W^{1,\infty}(-1,1)$ such that
  \begin{align*}
  &G(v)u_k=g_k(v) \qquad\forall\,v\in \mathcal B_1\,,\quad\forall\,k\in\enne\,,\\
  &g_k\sqrt{F''},\, g_k\sqrt{M''}\in L^\infty(-1,1) \quad\forall\,k\in\enne\,,\\
  &L_G:=\sum_{k=0}^\infty\left(\norm{g_k}_{W^{1,\infty}(-1,1)}^2 +
  \norm{g_k\sqrt{F''}}_{L^\infty(-1,1)}^2 + \norm{g_k\sqrt{M''}}_{L^\infty(-1,1)}^2\right)<+\infty\,.
  \end{align*}
  \item[{\bf D4}] The initial datum is nonrandom and satisfies
  \[
  \varphi_0\in V_1\,,\qquad |\varphi_0|<1 \text{ a.e.~in } \OO\,,\qquad
  F(\varphi_0)\in L^1(\OO)\,, \qquad M(\varphi_0)\in L^1(\OO)\,.
  \]
\end{itemize}
Note that under assumption {\bf D1} the irregular component 
of the potential $F$ is the convex part $F_1$, which may 
explode at $\pm1$. In condition {\bf D2} we assume 
on the other hand that the degeneracy of the mobility 
can only occur at $\pm1$, and compensates 
the eventual blow up of $F''$ at $\pm1$.
This is expected from the point of view
of application to phase-field modelling
(see the Remark~\ref{ex:log} below).
Finally, the additional summability condition in {\bf D3} is a generalization
of the classical compatibility condition between $m$ and $F''$
to the stochastic framework. This can be interpreted 
as a compensation of the blow up of $F''$ and $M''$ in $\pm1$ also
by the component functions $(g_k)_k$. Again,
this condition is satisfied in several 
physically relevant scenarios (see Remark~\ref{ex:log} below).

\begin{remark}[Logarithmic potential]
\label{ex:log}
Let us show now that the assumptions {\bf D1}--{\bf D2} allow for
the physically relevant case of degenerate mobility and 
logarithmic potential given by the natural choices
$m_{pol}$ and $F_{log}$ defined in \eqref{m_pol} and \eqref{pot_log},
respectively.
Indeed, assumption {\bf D1} holds with obvious choice of $F_1$ and $F_2$.
Moreover, an elementary computation yields
\[
  F_{log}''(r)=\frac\theta{1-r^2} - \theta_0\,, \quad r\in(-1,1)\,,
\]
so that $m_{pol}F_{log}''\in C^0([-1,1])$ and also condition {\bf D2} is satisfied.

\noindent With mobility $m_{pol}$ and potential $F_{log}$, 
a sufficient condition for 
assumption {\bf D3} is that 
\beq
\label{ex:noise}
  (g_k)_{k\in\enne}\subset W^{1,\infty}(-1,1)\,, \qquad
  g_k(-1)=g_k(1)=0 \quad\forall\,k\in\enne\,, \quad
  \sum_{k=0}^\infty\norm{g_k'}_{L^\infty(-1,1)}^2<+\infty\,,
\eeq
meaning essentially
that the components $(g_k)_k$ are Lipschitz-continuous
and vanish at the extremal points. Let us show that under \eqref{ex:noise}
also {\bf D3} is satisfied. Indeed,
for every $r\in(-1,1)$ and $k\in\enne$ one has
\begin{align*}
\left|g_k(r)\sqrt{F_{log}''(r)}\right|^2&=
\left|\theta\frac{g_k(r)}{1-r^2} - \theta_0g_k(r)\right|^2
\leq 2\theta\frac{g_k^2(r)}{1-r^2} + 2\theta_0 g_k^2(r)\\
&=2\theta\frac{|g_k(r)-g_k(-1)||g_k(r)-g_k(1)|}{|1-r||1+r|} +2\theta_0 |g_k(r)-g_k(1)|^2\\
&\leq2\theta\norm{g_k'}_{L^\infty(-1,1)}^2\frac{|1+r||1-r|}{|1+r||1-r|} + 
2\theta_0\norm{g_k'}^2_{L^\infty(-1,1)}|r-1|^2\\
&\leq 2\norm{g_k'}^2_{L^\infty(-1,1)}\left(\theta + 4\theta_0\right)\,,
\end{align*}
so that $g_k\sqrt{F_{log}''}\in L^\infty(-1,1)$ for all $k\in\enne$,
and
\[
\sum_{k=0}^\infty\norm{g_k\sqrt{F_{log}''}}_{L^\infty(-1,1)}^2\leq
2(\theta+4\theta_0)\sum_{k=0}^\infty\norm{g_k'}_{L^\infty(-1,1)}^2<+\infty\,.
\]
The computations for the terms $g_k\sqrt{M''_{pol}}$ are entirely analogous, 
and {\bf D3} follows.
\end{remark}

\begin{remark}[Double-obstacle potential]
  Note that choosing $F_1=0$ and $F_2(r):=1-r^2$, $r\in[-1,1]$,
  one recovers exactly the double-obstacle potential $F_{ob}$ defined in \eqref{F_ob}.
  Choosing also the degenerate mobility $m_{pol}$
  and $G$ as in \eqref{m_pol} and \eqref{ex:noise}, it is not difficult to
  check that {\bf D1}--{\bf D4} are satisfied.
\end{remark}

We give now the definition of martingale solution
in the case of degenerate mobility and irregular potential.
The main idea is to formally substitute equation \eqref{eq2}
in the variational formulation of the problem in order
to remove the dependence on the variable $\mu$.
The advantage of the degeneracy of the mobility 
is that the resulting variational formulation makes sense 
thanks to assumption {\bf D2}. 

\begin{defin}\label{def:sol-D}
  Assume {\bf D1}--{\bf D4}.
  A martingale solution to the problem \eqref{eq1}--\eqref{eq4} is a sextuple
  $(\hat\Omega, \hat\cF, (\hat\cF_t)_{t\in[0,T]}, \hat\P, \hat W, \hat\varphi)$ such that:
  \begin{itemize}
  \item $(\hat\Omega, \hat\cF, (\hat\cF_t)_{t\in[0,T]}, \hat\P)$ is a filtered probability 
  space satisfying the usual conditions;
  \item $\hat W$ is a $U$-valued cylindrical Wiener process on 
  $(\hat\Omega, \hat\cF, (\hat\cF_t)_{t\in[0,T]}, \hat\P)$;
  \item $\hat\varphi\in C^0_w([0,T]; L^2(\hat\Omega; H))\cap
  L^0(\hat\Omega; L^2(0,T; V_2))$ is progressively measurable;
  \item $|\hat\varphi|\leq1$ almost everywhere in $\hat\Omega\times Q$;
  \item it holds that 
  \beq\label{var-for-D}\begin{split}
  \int_\OO\hat\varphi(t,x)v(x)\,\d x &+
  \int_{Q_t}\Delta\hat\varphi(s,x)\div\left[m(\hat\varphi(s,x))\nabla v(x)\right]\,\d x\,\d s\\
  &+\int_{Q_t}m(\hat\varphi(s,x))F''(\hat\varphi(s,x))\nabla\hat\varphi(s,x)\cdot\nabla v(x)\,\d x\,\d s \\
  &=\int_\OO\varphi_0(x)v(x)\,\d x + 
  \int_\OO\left(\int_0^tG(\hat\varphi(s))\,\d \hat W(s)\right)\!\!(x)v(x)\,\d x
  \end{split}\eeq
  for every $v\in V_2\cap W^{1,d}(\OO)$, $\hat\P$-almost surely, for every $t\in[0,T]$.
  \end{itemize}
\end{defin}

\begin{remark}
  As we have anticipated, the variational formulation \eqref{var-for-D} 
  in Definition~\ref{def:sol-D}
  is formally obtained substituting the definition \eqref{eq2} of chemical potential
  in equation \eqref{eq1} and integrating by parts.
  The reason why we do so is that the chemical potential $\mu$
  does not inherit enough regularity in the degenerate case.
  The main advantage of such substitution is that 
  all the terms in \eqref{var-for-D} are still well-defined.
  Indeed, by {\bf D2} we have that $mF''\in C^0([-1,1])$, so that 
  the third term on the left-hand side makes sense. Moreover, 
  for every $v\in V_2\cap W^{1,d}(\OO)$ we have that 
  $\div(m(\hat\varphi)\nabla v)=m'(\hat\varphi)\nabla\hat\varphi\cdot\nabla v + m(\hat\varphi)\Delta v$,
  where $m(\hat\varphi), m'(\hat\varphi)\in L^\infty(Q)$ by {\bf D2}.
  Also, since $V_1\embed L^{\frac{2d}{d-2}}(\OO)$,
  noting that $\frac{d-2}{2d} + \frac1d = \frac12$ by the H\"older inequality 
  we have
  $\div(m(\varphi')\nabla v)\in H$, so that also the second term on the left-hand side of 
  \eqref{var-for-D} makes sense.
\end{remark}

We are now ready to state the second main result of the paper, 
ensuring existence of martingale solutions in the case of 
degenerate mobility and irregular potential. Both the cases of logarithmic 
and double-obstacle potential are covered.

\begin{thm}
  \label{thm2}
  Assume conditions {\bf D1}--{\bf D4}. Then,
  there exists a martingale solution
  \[
  \left(\hat\Omega, \hat\cF, (\hat\cF_t)_{t\in[0,T]}, \hat\P, \hat W, \hat\varphi\right)
  \]
  to the problem \eqref{eq1}--\eqref{eq4} in the sense of Definition~\ref{def:sol-D}
  such that
  \begin{align*}
  &\hat\varphi \in C^0_w([0,T]; L^2(\hat \Omega; V_1))\cap L^2(\hat\Omega; L^2(0,T; V_2))\,,\\
  &F(\hat\varphi),\,M(\hat\varphi)\in L^\infty(0,T; L^1(\Omega\times\OO))\,.
  \end{align*}
  In particular, if
  \beq
  \label{sing_lim}
  \lim_{|r|\to 1^-} F_1(r)=+\infty \qquad\text{or}\qquad
  \lim_{|r|\to 1^-} M(r)=+\infty\,,
  \eeq
  then
  \[
  |\hat\varphi(t)|<1 \quad\text{a.e.~in } \hat\Omega\times\OO \quad\forall\,t\in[0,T]\,.
  \]
\end{thm}
Let us stress that the last assertion of Theorem~\ref{thm2} 
ensures that under \eqref{sing_lim} the concentration $\hat\varphi$
is almost everywhere contained in the interior of the physically relevant domain,
meaning that the contact set $\{|\hat\varphi|=1\}$ has measure $0$, or better said that
\[
  |\{(\hat\omega,x)\in\hat\Omega\times \OO:|\hat\varphi(\hat\omega,t,x)|=1\}|=0 \quad\forall\,t\in[0,T]\,.
\]
Note that in general
the degeneracy of the mobility at $\pm1$
may prevent $M$ to blow up at $\pm1$,
and \eqref{sing_lim} is not always satisfied.
For example, 
an easy computation shows that for degenerate mobility $m_{pol}$
introduced in Remark~\ref{ex:log} we have
\[
  M_{pol}(r)=\frac12\left((1+r)\ln(1+r)+(1-r)\ln(1-r)\right)\,, \quad r\in(-1,1)\,.
\]
Hence, in such a case $M_{pol}$ is bounded in $(-1,1)$
and condition \eqref{sing_lim} is satisfied 
only if the potential $F$ blows up at $\pm1$.
On the other hand, 
in case of polynomial mobility $m_\alpha$
with $\alpha\geq2$,
condition \eqref{sing_lim} is always satisfied, 
irrespectively of the potential $F$.
If \eqref{sing_lim} is not satisfied, then one can only infer
that $|\hat\varphi|\leq1$ almost everywhere, as it is natural to expect.

%%%%%%%%%%%%%%%%%%%%%%%%%%%%%%%%%%%%%%%%%%%%%%%%%%%%%

\section{Positive mobility and regular potential}
\label{sec:ND}
This section is devoted to the proof of Theorem~\ref{thm1}.
The main idea is to perform
two separate approximations on the problem.
The first one depends on a 
parameter $\lambda>0$,
and is obtained replacing
the nonlinearity $F$ with 
its Yosida approximation.
The second one
depends on the parameter $n\in\enne$ and
is a Faedo-Galerkin
finite-dimensional approximation.
Uniform estimates are proved
first uniformly in $n$, when
$\lambda$ is fixed,
and a passage to the limit 
as $n\to\infty$ yields
existence of approximated solutions
for $\lambda>0$ fixed.
Secondly, further uniform estimates
are proved uniformly in $\lambda$
and a passage to the limit as $\lambda\to0$
gives existence of solutions to the original problem.

\subsection{The approximation}
First of all, since $F''\geq-C_F$, the function
$\gamma:\erre\to\erre$
defined as $\gamma(r):=F'(r) + C_Fr$, $r\in\erre$, is nondecreasing and continuous:
hence, $\gamma$ can be identified with a maximal monotone graph
in $\erre\times\erre$ and satisfies $\gamma(0)=0$.
It makes sense then to introduce
the Yosida approximation $\gamma_\lambda:\erre\to\erre$
of $\gamma$ for any $\lambda>0$ and define
$\hat\gamma_\lambda:\erre\to\mathopen[0,+\infty\mathclose)$
as $\hat\gamma_\lambda(r):=\int_0^r\gamma_\lambda(s)\,\d s$, $r\in\erre$.
With this notation, we introduce the approximated potential as
\[
  F_\lambda:\erre\to\mathopen[0,+\infty\mathclose)\,, \qquad
  F_\lambda(r):=F(0) + \hat\gamma_\lambda(r) - \frac{C_F}2r^2\,, \quad r\in\erre\,.
\]
Let us recall that from the general theory on monotone analysis \cite[Ch.~2]{barbu-monot}
we know that the Yosida approximation $\gamma_\lambda$ is 
$\frac1\lambda$-Lipschitz-continuous, hence also linearly bounded. 
Consequently, by definition of $\hat\gamma_\lambda$ there exists a constant 
$c_\lambda>0$ such that $\hat\gamma(r)\leq c_\lambda(1+|r|^2)$
for all $r\in\erre$. Also, noting that $F_\lambda'(r)=\gamma_\lambda(r)-C_Fr$
for all $r\in\erre$ by definition, 
this readily ensures that also $F_\lambda'$ is $\frac1\lambda$-Lipschitz-continuous
and that, possibly renominating $c_\lambda$, it holds
\beq
  \label{Flam_quad}
  |F_\lambda(r)|\leq c_\lambda(1+|r|^2) \qquad\forall\,r\in\erre\,.
\eeq

Secondly, let $(e_j)_{j\in\enne_+}\subset V_2$ and $(\alpha_j)_{j\in\enne_+}$ 
be the sequences of eigenfunctions and eigenvalues of the 
negative Laplace operator with homogeneous Neumann conditions on $\OO$, respectively, i.e.
\[
  \begin{cases}
  -\Delta e_j = \alpha_j e_j \quad&\text{in } \OO\,,\\
  {\bf n}\cdot \nabla e_j = 0 \quad&\text{in } \partial\OO\,,
  \end{cases}
  \qquad j\in\enne_+\,.
\]
Then, possibly using a renormalization procedure, 
we can suppose that $(e_j)_j$ is a complete orthonormal system of $H$
and an orthogonal system in $V_1$.
For every $n\in\enne_+$, we define the finite
dimensional space $H_n:=\operatorname{span}\{e_1,\ldots,e_n\}\subset V_2$,
endowed with the $\norm{\cdot}_H$-norm.

We define the approximated operator $G_{n}:H_n\to\cL^2(U,H_n)$ as
\[
  G_{n}(v)u_k:=\sum_{j=1}^n(G(v)u_k, e_j)_H e_j\,, \qquad v\in H_n\,, \quad k\in\enne\,.
\]
One can check that $G_{n}$ is well-defined: indeed, for every $v\in H_n$ 
and every $n\in\enne_+$, thanks to assumption {\bf ND3} we have 
\begin{align*}
  \sum_{k=0}^\infty\norm{G_{n}(v)u_k}_H^2 &= 
  \sum_{k=0}^\infty\sum_{j=1}^n|(G(v)u_k, e_j)_H|^2\leq
  \sum_{k=0}^\infty\sum_{j=1}^\infty|(G(v)u_k, e_j)_H|^2\\ 
  &=\sum_{k=0}^\infty\norm{G(v)u_k}_H^2= \norm{G(v)}^2_{\cL^2(U,H)}\,,
\end{align*}
so that $G_{n}(v)\in\cL^2(U,H_n)$ for all $v\in H_n$ and
\beq
  \label{est_Gn}
  \norm{G_{n}(v)}_{\cL^2(U,H)}\leq
  \norm{G(v)}_{\cL^2(U,H)} \qquad\forall\,v\in H,
  \quad\forall\,n\in\enne_+\,.
\eeq
A similar computations shows also that $G_{n}$ is Lipschitz-continuous
from $H_n$ to $\cL^2(U,H_n)$.

Similarly, we define the approximated initial value
\[
  \varphi_0^{n}:=\sum_{j=1}^n(\varphi_0, e_j)_He_j\,.
\]

Finally, 
for every $n\in\enne_+$ let $m_n:=\rho_n*m$
where $(\rho_n)$ is a standard sequence of mollifiers.
In particular, we have that
\[
  (m_n)_n\subset W^{1,\infty}(\erre)\,, \qquad
  m_*\leq m_n(r)\leq m^*\quad\forall\,r\in\erre\,,\qquad
  m_n\to m \quad\text{in } C^0([a,b])\quad\,\forall a<b\,.
\]

We consider the approximated problem
\begin{align}
  \label{eq1_app}
  \d\varphi_{\lambda,n} - \div\left(m_n(\varphi_{\lambda,n})\nabla\mu_{\lambda,n}\right)\,\d t
  =G_{n}(\varphi_{\lambda,n})\,\d W
  \qquad&\text{in } (0,T)\times\OO\,,\\
  \label{eq2_app}
  \mu_{\lambda,n}=-\Delta\varphi_{\lambda,n} + F_{\lambda}'(\varphi_{\lambda,n})
  \qquad&\text{in } (0,T)\times\OO\,,\\
  \label{eq3_app}
  {\bf n}\cdot\nabla\varphi_{\lambda,n} = 
  {\bf n}\cdot m(\varphi_{\lambda,n})\nabla\mu_{\lambda,n} = 0
  \qquad&\text{in } (0,T)\times\partial\OO\,,\\
  \label{eq4_app}
  \varphi_{\lambda,n}(0)=\varphi_0^{n}
  \qquad&\text{in } \OO\,.
\end{align}
Let us fix now$\lambda>0$ and $n\in\enne_+$:
we look for a solution $(\varphi_{\lambda,n},\mu_{\lambda,n})$ 
to \eqref{eq1_app}--\eqref{eq4_app} in the form
\[
  \varphi_{\lambda,n}=\sum_{j=1}^n a_j^{\lambda,n} e_j\,, \qquad
  \mu_{\lambda,n}=\sum_{j=1}^n b_j^{\lambda,n} e_j\,,
\]
for some processes 
\[
a^{\lambda,n}:=(a_1^{\lambda,n},\ldots, a_n^{\lambda,n}):\Omega\times[0,T]\to\erre^n\,, \qquad
b^{\lambda,n}:=(b_1^{\lambda,n},\ldots,b_n^{\lambda,n}):\Omega\times[0,T]\to\erre^n\,.
\]
Plugging in the {\em ansatz} on $\varphi_{\lambda,n}$ and $\mu_{\lambda,n}$ 
in \eqref{eq1_app}--\eqref{eq4_app} and
taking any arbitrary $e_i$, $i=1,\ldots,n$, as test functions, we immediately see that 
the variational formulation of \eqref{eq1_app}--\eqref{eq4_app} is given by 
\begin{align*}
  &\int_\OO\varphi_{\lambda,n}(t,x)e_i(x)\,\d x 
  + \int_{Q_t}m_n(\varphi_{\lambda,n}(s,x))\nabla\mu_{\lambda,n}(s,x)\cdot\nabla e_i(x)\,\d x\,\d s\\
  &\qquad=\int_\OO\varphi_0^{n}(x)e_i(x)\,\d x + 
  \int_\OO\left(\int_0^tG_{n}(\varphi_{\lambda,n}(s))\,\d W(s)\right)\!\!(x)e_i(x)\,\d x \qquad
  \forall\,i=1,\ldots,n\,,
\end{align*}
and
\[
  \int_\OO\mu_{\lambda,n}(t,x)e_i(x)\,\d x 
  = \int_\OO\nabla\varphi_{\lambda,n}(t,x)\cdot\nabla e_i(x)\,\d x
  +\int_\OO F'_\lambda(\varphi_{\lambda,n}(t,x))e_i(x)\,\d x
  \qquad\forall\,i=1,\ldots,n\,,
\]
for every $t\in[0,T]$, $\P$-almost surely.
Using the orthogonality properties of $(e_j)_j$,
we deduce then that 
$(\varphi_{\lambda,n}, \mu_{\lambda,n})$ satisfy \eqref{eq1_app}--\eqref{eq4_app} 
if and only if the vectors $(a^{\lambda,n}, b^{\lambda,n})$ satisfy the SDEs
\begin{align*}
  &\d a^{\lambda,n}_i + \sum_{j=1}^nb_j^{\lambda,n}
  \int_\OO m_n\left(\sum_{l=1}^n a^{\lambda,n}_l e_l(x)\right)\nabla e_j(x)\cdot\nabla e_i(x)\,\d x\,\d t
  =\left(G_{n}\left(\sum_{l=1}^n a^{\lambda,n}_l e_l\right)\,\d W, e_i\right)_H\,,\\
  &b^{\lambda,n}_i=\alpha_i a^{\lambda,n}_i 
  + \int_\OO F'_\lambda\left(\sum_{l=1}^n a^{\lambda,n}_l e_l(x)\right)e_i(x)\,\d x\,,\\
  &a^{\lambda,n}_i(0)=\left(\varphi_0, e_i\right)_H\,,
\end{align*}
for every $i=1,\ldots,n$, where
the stochastic integral on the right-hand side must be interpreted
for any $i=1,\ldots,n$ as $G_{i}^{\lambda,n}\,\d W$ where
\[
  G^{\lambda,n}_{i}:H_n\to \cL^2(U,\erre)\,, \qquad
  G^{\lambda,n}_{i}u_k:=
  \left(G_{n}\left(\sum_{l=1}^n a^{\lambda,n}_l e_l\right)u_k, e_i\right)_H\,, \quad k\in\enne\,.
\]
Since the functions $m_n$, $F_\lambda'$, 
and the operator $G_{n}$ are Lipschitz-continuous,
the system can be seen as an abstract evolution equation
on the Hilbert space $H_n$ with Lipschitz--continuous nonlinearities:
hence, the classical theory \cite[\S~7.1]{dapratozab} applies
and we can find a unique solution
\[
  a^{\lambda,n}, b^{\lambda,n} \in L^\ell\left(\Omega; C^0([0,T]; \erre^n)\right)
  \qquad\forall\,\ell\in\mathopen[2,+\infty\mathclose)\,.
\]
We deduce that for every $n\in\enne$ 
the approximated system \eqref{eq1_app}--\eqref{eq4_app}
admits a unique solution
\[
  \varphi_{\lambda,n}, \mu_{\lambda,n} \in L^\ell\left(\Omega; C^0([0,T]; H_n)\right)
  \qquad\forall\,\ell\in\mathopen[2,+\infty\mathclose)\,.
\]

\subsection{Uniform estimates in $n$, with $\lambda$ fixed}
We show now that the approximated solution $(\varphi_{\lambda,n}, \mu_{\lambda,n})$
satisfy some energy estimates, independently of $n$,
with $\lambda>0$ being fixed.

First of all, integrating \eqref{eq1_app} on $\OO$ and using It\^o's formula yields
\begin{align*}
  \frac12|(\varphi_{\lambda,n}(t))_\OO|^2 &= \frac12|(\varphi_0^{n})_\OO|^2
  +\int_0^t(\varphi_{\lambda,n}(s))_\OO(G_{n}(\varphi_{\lambda,n}(s))\d W(s))_\OO\\
  &+\frac12\int_0^t\norm{(G_{n}(\varphi_{\lambda,n}(s))_\OO}_{\cL^2(U,\erre)}^2\,\d s\,.
\end{align*}
Taking supremum in time, power $\ell/2$ and expectations, thanks to the 
Burkholder-Davis-Gundy inequality we have
\begin{align*}
  \E\sup_{s\in[0,t]}|(\varphi_{\lambda,n}(s))_\OO|^\ell
  &\lesssim \E|(\varphi_0^{n})_\OO|^\ell
  +\E\left(\int_0^t\sum_{k=0}^\infty|(G_{n}(\varphi_{\lambda,n}(s))u_k)_\OO|^2\,\d s\right)^{\ell/2}\\
  &+\E\left(\int_0^t|(\varphi_{\lambda,n}(s))_\OO|^2
  \sum_{k=0}^\infty|(G_{n}(\varphi_{\lambda,n}(s))u_k)_\OO|^2\,\d s\right)^{\ell/4}\,.
\end{align*}
Note that $|(\varphi_0^{n})_\OO|=|\OO|^{-1}|\int_\OO\varphi_0^n|\leq
|\OO|^{-1}\norm{\varphi_0^n}_{L^1(\OO)}
\leq|\OO|^{-1/2}
\|\varphi_0^{n}\|_H\leq|\OO|^{-1/2}\norm{\varphi_0}_H$
and 
\begin{align*}
\sum_{k=0}^\infty|(G_{n}(\varphi_{\lambda,n})u_k)_\OO|^2&\leq
|\OO|^{-1}\sum_{k=0}^\infty\norm{G_{n}(\varphi_{\lambda,n})u_k}_H^2\leq
|\OO|^{-1}\sum_{k=0}^\infty\norm{G(\varphi_{\lambda,n})u_k}_H^2\\
&=|\OO|^{-1}\sum_{k=0}^\infty\norm{g_k(\varphi_{\lambda,n})}_H^2
\leq \sum_{k=0}^\infty\norm{g_k}_{L^\infty(\erre)}^2\leq C_G\,.
\end{align*}
Hence, by the Young inequality we infer that 
there exists $c>0$, independent of $\lambda$ and $n$, such that 
\beq
  \label{est_mean}
  \norm{(\varphi_{\lambda,n})_\OO}_{L^\ell(\Omega; C^0([0,T]))}\leq c\,.
\eeq

We want now to write It\^o's formula for the free energy functional
\beq\label{energy}
  \mathcal E_\lambda(v):=\frac12\int_\OO|\nabla v|^2 + \int_\OO F_\lambda(v)\,, \qquad v\in H_n\,.
\eeq
To this end, note that 
since $H_n\embed V_2$ and $F_\lambda'$ is Lipschitz-continuous,
$\mathcal E_\lambda:H_n\to\mathopen[0,+\infty\mathclose)$
is Fr\'echet-differentiable with $D\mathcal E_\lambda: H_n\to H_n^*$
given by
\[
  D\mathcal E_\lambda(v)h=\int_\OO\nabla v(x)\cdot\nabla h(x)\,\d x + \int_\OO F_\lambda'(v(x))h(x)\,\d x\,,
  \quad v,h\in H_n\,.
\] 
Let us show now that also $D\mathcal E_\lambda$ is Fr\'echet-differentiable with 
$D^2\mathcal E_\lambda: H_n\to \cL(H_n,H_n^*)$ given by 
\[
  D^2\mathcal E_\lambda(v)[h,k]=\int_\OO\nabla h(x)\cdot\nabla k(x)\,\d x
  +\int_\OO F_\lambda''(v(x))h(x)k(x)\,\d x\,, \quad v,h,k\in H_n\,.
\]
Indeed, for every $v,h,k\in H_n$
we have that 
\begin{align*}
  &\left|D\mathcal E_\lambda(v+k)h - D\mathcal E_\lambda(v)h - 
  \int_\OO\nabla h(x)\cdot\nabla k(x)\,\d x
  -\int_\OO F_\lambda''(v(x))h(x)k(x)\,\d x\right|\\
  &=\left|\int_\OO F_\lambda'(v(x)+k(x))h(x)\,\d x - \int_\OO F_\lambda'(v(x))h(x)\,\d x
  -\int_\OO F_\lambda''(v(x))h(x)k(x)\,\d x\right|\\
  &=\left|\int_0^1\int_\OO \left(F_\lambda''(v(x) + \tau k(x)) - 
  F_\lambda''(v(x))\right)h(x)k(x)\,\d x\, \d \tau\right|\,.
\end{align*}
Now, since we have the continuous inclusion
$H_n\embed L^{\infty}(\OO)$,
by H\"older inequality we infer that 
\begin{align*}
  &\left|D\mathcal E_\lambda(v+k)h - D\mathcal E_\lambda(v)h - 
  \int_\OO\nabla h(x)\cdot\nabla k(x)\,\d x
  -\int_\OO F_\lambda''(v(x))h(x)k(x)\,\d x\right|\\
  &\leq\norm{h}_{L^{\infty}(\OO)}\norm{k}_{L^{\infty}(\OO)}
  \int_0^1\norm{F_\lambda''(v+\tau k) - F_\lambda''(v)}_{L^1(\OO)}\,\d \tau\\
  &\leq c_n^2\norm{h}_{H_n}\norm{k}_{H_n}
  \int_0^1\norm{F_\lambda''(v+\tau k) - F_\lambda''(v)}_{L^{1}(\OO)}\,\d \tau\,,
\end{align*}
where $c_n>0$ is the norm of the inclusion $H_n\embed L^{\infty}(\OO)$.
Since $F_\lambda''$ is continuous and bounded, the third factor on the right-hand side
converges to $0$ if $k\to 0$ in $H_n$ by the dominated convergence theorem, hence
\begin{align*}
  &\sup_{\norm{h}_{H_n}\leq 1}
  \left|D\mathcal E_\lambda(v+k)h - D\mathcal E_\lambda(v)h - 
  \int_\OO\nabla h(x)\cdot\nabla k(x)\,\d x
  -\int_\OO F_\lambda''(v(x))h(x)k(x)\,\d x\right| \\
  &=o\left(\norm{k}_{H_n}\right)
  \quad\text{as } k\to 0 \text{ in } H_n\,.
\end{align*}
This shows that $D\mathcal E_\lambda$ is indeed Fr\'echet-differentiable with derivative 
$D^2\mathcal E_\lambda$ given as above.
Furthermore, the derivatives $D\mathcal E_\lambda$ and $D^2 \mathcal E_\lambda$
are continuous and bounded on bounded subsets of $H_n$, as it follows
from the Lipschitz-continuity of $F_\lambda'$ and the continuity and boundedness of $F_\lambda''$.

We can then apply It\^o's formula to $\mathcal E_\lambda(\varphi_{\lambda,n})$ 
in the classical version of \cite{dapratozab}.
To this end, note that by \eqref{eq2_app} we have 
$D\mathcal E_\lambda(\varphi_{\lambda,n})=\mu_{\lambda,n}$: we obtain then
\beq\label{ito_app}\begin{split}
  &\mathcal E_\lambda(\varphi_{\lambda,n}(t))
  +\int_{Q_t}m_n(\varphi_{\lambda,n}(s,t))|\nabla\mu_{\lambda,n}(s,x)|^2\,\d x\,\d s\\
  &=\mathcal E_\lambda(\varphi_0^{n}) 
  +\int_0^t\left(\mu_{\lambda,n}(s), G_{n}(\varphi_{\lambda,n}(s))\,\d W(s)\right)_H\\
  &+\frac12\int_0^t\sum_{k=0}^\infty\int_\OO
  \left[|\nabla G_{n}(\varphi_{\lambda,n}(s))u_k|^2(x)
  +F_\lambda''(\varphi_{\lambda,n}(s,x))|G_{n}(\varphi_{\lambda,n}(s))u_k|^2(x)\right]\,\d x\,\d s\,.
\end{split}\eeq
Taking power $\ell/2$ at both sides, supremum in time and then expectations yields,
recalling the definition \eqref{energy} and the assumption {\bf ND2},
\begin{align*}
   &\frac12\E\sup_{s\in[0,t]}\norm{\nabla\varphi_{\lambda,n}(s)}_H^{\ell} +
   \E\sup_{s\in[0,t]}\norm{F_\lambda(\varphi_{\lambda,n}(s))}_{L^1(\OO)}^{\ell/2}
   +m_*^{\ell/2}\E\norm{\nabla\mu_{\lambda,n}}_{L^2(0,t; H)}^{\ell}\\
   &\leq c_\ell\left[
   \norm{\nabla\varphi_0^{n}}_H^{\ell} 
   +\norm{F_\lambda(\varphi_0^{n})}_{L^1(\OO)}^{\ell/2}
   +\E\sup_{s\in[0,t]}\left|\int_0^s\left(\mu_{\lambda,n}(r), 
   G_{n}(\varphi_{\lambda,n}(r))\,\d W(r)\right)_H\right|^{\ell/2}\right.\\
   &\left.+\E\norm{G_{n}(\varphi_{\lambda,n})}^{\ell}_{L^2(0,t; \cL^2(U,V_1))}
   +\E\left|\int_0^t\sum_{k=0}^\infty\int_\OO
  |F_\lambda''(\varphi_n(s,x))||G_{n}(\varphi_{\lambda,n}(s))u_k|^2(x)\,\d x\,\d s\right|^{\ell/2}
  \right]
\end{align*}
for every $t\in[0,T]$, $\P$-almost surely, for a certain constant $c_\ell>0$
depending only on $\ell$.
Let us estimate the terms on the right-hand side separately.
First of all, from the definition of the approximate initial value
$\varphi_0^n$, since $\varphi\in V_1$ we have 
$\|\nabla\varphi_0^{n}\|_H\leq\norm{\nabla\varphi_0}_H$.
Let us focus on the stochastic integral.
By the Burkholder-Davis-Gundy inequality and
the estimate \eqref{est_Gn}, 
we infer that 
\begin{align*}
  \E\sup_{s\in[0,t]}\left|\int_0^s\left(\mu_{\lambda,n}(r), 
  G_{n}(\varphi_n(r))\,\d W(r)\right)_H\right|^{\ell/2}
  &\leq c_\ell
  \E\left(\int_0^t\norm{\mu_{\lambda,n}(s)}_H^2
  \norm{G_{n}(\varphi_{\lambda,n})}_{\cL^2(U,H)}^2\,\d s\right)^{\ell/4}\\
  &\leq c_\ell
  \E\left(\int_0^t\norm{\mu_{\lambda,n}(s)}_H^2\norm{G(\varphi_{\lambda,n})}_{\cL^2(U,H)}^2\,\d s\right)^{\ell/4}
\end{align*}
for a certain constant $c_\ell>0$ independent of $n$ and $\lambda$.
Thanks to assumption {\bf ND3} we have
\[
  \norm{G(\varphi_{\lambda,n})}_{\cL^2(U,H)}^2=\sum_{k=0}^\infty\norm{g_k(\varphi_{\lambda,n})}_H^2
  \leq|\OO|\sum_{k=0}^\infty\norm{g_k}_{L^\infty(\erre)}^2\leq |O|C_G\,,
\]
so that, consequently, we deduce that 
\begin{align*}
  \E\sup_{s\in[0,t]}\left|\int_0^s\left(\mu_{\lambda,n}(r), 
  G_{n}(\varphi_{\lambda,n}(r))\,\d W(r)\right)_H\right|^{\ell/2}
  &\leq c_\ell |O|^{\ell/4}C_G^{\ell/4}
  \E\norm{\mu_{\lambda,n}}_{L^2(0,t; H)}^{\ell/2}\,.
\end{align*}
Summing and subtracting $(\mu_{\lambda,n})_\OO$ on the right-hand side,
using the Poincar\'e-Wirtinger inequality and the Young inequality we deduce that 
\begin{align*}
  &\E\sup_{s\in[0,t]}\left|\int_0^s\left(\mu_{\lambda,n}(r), 
  G_{n}(\varphi_{\lambda,n}(r))\,\d W(r)\right)_H\right|^{\ell/2}\\
  &\leq \frac{m_*^{\ell/2}}{2} \E\norm{\nabla\mu_{\lambda,n}}_{L^2(0,t; H)}^{\ell} 
  + c\E\norm{(\mu_{\lambda,n})_\OO}_{L^2(0,t)}^{\ell/2} + c\,,
\end{align*}
where $c=c(m_*, c_\ell, |\OO|, \ell, C_G)>0$ is an arbitrarily large constant
independent of $n$ and $\lambda$.
Let us focus now on the trace terms in It\^o's formula.
Since $G(\varphi_{\lambda,n})$ takes values in $\cL^2(U,V_1)$, we have
\begin{align*}
  \norm{G_{n}(\varphi_{\lambda,n})}_{\cL^2(U,V_1)}^2
  &\leq\norm{G(\varphi_{\lambda,n})}_{\cL^2(U,V_1)}^2
  =\sum_{k=0}^\infty\norm{g_k(\varphi_{\lambda,n})}_{V_1}^2\\
  &\leq\left(|\OO| + \norm{\nabla\varphi_{\lambda,n}}_H^2\right) 
  \sum_{k=0}^\infty\norm{g_k}^2_{W^{1,\infty}(\erre)}\,,
\end{align*}
so that {\bf ND3} yields
\[
  \E\norm{G_{n}(\varphi_n)}^\ell_{L^2(0,t; H)} 
  \leq c\left(1 + \E\norm{\nabla\varphi_{\lambda,n}}^\ell_{L^2(0,t; H)}\right)
\]
for a certain constant $c>0$ independent of $n$ and $\lambda$. 
Taking the mean of \eqref{eq2_app} we get 
\[
  (\mu_{\lambda,n})_\OO = (F_\lambda'(\varphi_{\lambda,n}))_\OO\leq
  c\norm{F_\lambda'(\varphi_{\lambda,n})}_{L^1(\OO)}\,,
\]
so that from assumption {\bf ND1} it follows that
\[
  |(\mu_{\lambda,n})_\OO|\leq c\left(1 + \norm{F_\lambda(\varphi_{\lambda,n})}_{L^1(\OO)}\right)
\]
for a certain $c>0$ independent of $\lambda$ and $n$:
we infer then that, possibly updating the value of $c$,
\beq\label{ito_aux}\begin{split}
   &\E\sup_{s\in[0,t]}\norm{\nabla\varphi_{\lambda,n}(s)}_H^{\ell} +
   \E\sup_{s\in[0,t]}\norm{F_\lambda(\varphi_{\lambda,n}(s))}_{L^1(\OO)}^{\ell/2}\\
   &\qquad+\E\sup_{s\in[0,t]}|(\mu_{\lambda,n}(s))_\OO|^{\ell/2}
   +\E\norm{\nabla\mu_{\lambda,n}}_{L^2(0,t; H)}^{\ell}\\
   &\leq c\left(1 + \E\norm{(\mu_{\lambda,n})_\OO}_{L^2(0,t)}^{\ell/2}
   + \E\norm{\nabla\varphi_{\lambda,n}}^\ell_{L^2(0,t; H)}\right)
   +\norm{F_\lambda(\varphi_0^{n})}_{L^1(\OO)}^{\ell/2}\\
   &\qquad+\frac12\E\left(\int_0^t\sum_{k=0}^\infty\int_\OO
  |F_\lambda''(\varphi_{\lambda,n}(s,x))||G_{n}(\varphi_{\lambda,n}(s))u_k|^2(x)\,\d x\,\d s\right)^{\ell/2}\,.
\end{split}\eeq
Let us estimate the two terms on the right-hand side.
Recall that here $\lambda>0$ is fixed.
First of all, since $F_\lambda$ is bounded by a quadratic function by \eqref{Flam_quad}
and $(\varphi_0^n)_n$ is bounded in $H$ thanks to the properties of
the orthogonal projection on $H_n$, 
we have that 
\[
\norm{F_\lambda(\varphi_0^{n})}_{L^1(\OO)}\leq 
c_\lambda\left(1+\norm{\varphi_0^n}_H^2\right)\leq 
c_\lambda\left(1+\norm{\varphi_0}_H^2\right)
 \qquad\forall\,n\in\enne
\]
for a certain $c_\lambda>0$ independent of $n$.
Secondly, note that 
since $|F''_\lambda|\leq c_\lambda$ for a certain $c_\lambda>0$
independent of $n$, by the same computations as above we have
\[
  \sum_{k=0}^\infty\int_\OO
  |F_\lambda''(\varphi_{\lambda,n}(s,x))||G_{n}(\varphi_n(s))u_k|^2(x)\,\d x
  \leq c_\lambda\norm{G(\varphi_{\lambda,n})}_{\cL^2(U,H)}^2
  \leq c_\lambda |\OO| C_G\,.
\]

Putting this information together
we deduce from \eqref{ito_aux} that
there exists a positive constant $c_\lambda$, independent of $n$,
such that 
\begin{align*}
   &\E\!\sup_{s\in[0,t]}\norm{\nabla\varphi_{\lambda,n}(s)}_H^{\ell} +
   \E\!\sup_{s\in[0,t]}\norm{F_\lambda(\varphi_{\lambda,n}(s))}_{L^1(\OO)}^{\ell/2}+
   \E\!\sup_{s\in[0,t]}|(\mu_{\lambda,n}(s))_\OO|^{\ell/2}
   +\E\!\norm{\nabla\mu_{\lambda,n}}_{L^2(0,t; H)}^{\ell}\\
   &\leq c_\lambda\left(1 + \E\norm{(\mu_{\lambda,n})_\OO}_{L^2(0,t)}^{\ell/2}
   + \E\norm{\nabla\varphi_{\lambda,n}}^\ell_{L^2(0,t; H)}\right)\,.
\end{align*}
The Gronwall lemma and the estimate \eqref{est_mean}
yield then, after updating the constant $c_\lambda$,
\begin{align}
  \label{est1_n}
  \norm{\varphi_{\lambda,n}}_{L^\ell(\Omega; C^0([0,T]; V_1))}&\leq c_\lambda\,,\\
  \label{est2_n}
  \norm{\mu_{\lambda,n}}_{L^{\ell/2}(\Omega; L^2(0,T; V_1))} + 
  \norm{\nabla\mu_{\lambda,n}}_{L^\ell(\Omega; L^2(0,T; H))} &\leq c_\lambda\,.
\end{align}
Moreover, the computations performed above also imply that 
\[
  \norm{G_n(\varphi_{\lambda,n})}_{L^\infty(\Omega\times(0,T); \cL^2(U,H))
  \cap L^\ell(\Omega; L^\infty(0,T; \cL^2(U,V_1)))}\leq c_\lambda
\]
for a positive constant $c$ independent of $\lambda$ and $n$, hence also by
\cite[Lem.~2.1]{fland-gat}, for any $s\in(0,1/2)$,
\[
  \norm{\int_0^\cdot G_n(\varphi_{\lambda,n}(s))\,\d W(s)}_{L^\ell(\Omega; W^{s,\ell}(0,T; V_1))}\leq c_{\lambda,s}\,.
\]
By comparison in \eqref{eq2_app} we infer then that 
\beq\label{est3_n}
  \norm{\varphi_{\lambda,n}}_{L^\ell(\Omega; W^{\bar s,\ell}(0,T; V_1^*))}\leq c_\lambda\,,
\eeq
where $c_\lambda$ is independent of $n$, and 
$\bar s\in(1/\ell,1/2)$ is fixed.

\subsection{Passage to the limit as $n\to\infty$, with $\lambda$ fixed}
\label{ssec:lim_n}
We perform here the passage to the limit as $n\to\infty$, 
keeping $\lambda>0$ fixed.

Let us show that the sequence of laws of $(\varphi_{\lambda,n})_n$ is tight on $C^0([0,T]; H)$.
To this end, let us recall that, since $\bar s>1/\ell$,
by \cite[Cor.~5, p.~86]{simon}
we have the compact inclusion
\[
  L^\infty(0,T; V_1)\cap W^{\bar s,\ell}(0,T; V_1^*)\cembed C^0([0,T]; H)\,.
\]
Hence, for every $R>0$ the closed ball $B_R$ 
in $L^\infty(0,T; V_1)\cap W^{\bar s,\ell}(0,T; V_1^*)$ of radius $R$
is compact in $C^0([0,T]; H)$. Moreover, 
thanks to the Markov inequality and the estimates \eqref{est1_n} and \eqref{est3_n} 
we have
\begin{align*}
  \P\{\varphi_{\lambda,n}\in B_R^c\}&=
  \P\{\norm{\varphi_{\lambda,n}}_{L^\infty(0,T; V_1)\cap W^{\bar s,\ell}(0,T; V_1^*)}> R\}\\
  &\leq
  \frac1{R^{\ell}}\E\norm{\varphi_{\lambda,n}}_{L^\infty(0,T; V_1)\cap W^{\bar s,\ell}(0,T; V_1^*)}^\ell\leq
  \frac{c_\lambda^\ell}{R^\ell}\,,
\end{align*}
which yields
\[
  \lim_{R\to+\infty}\sup_{n\in\enne_+}\P\{\varphi_{\lambda,n}\in B_R^c\}=0\,,
\]
as required. Hence, the family of laws of $(\varphi_{\lambda,n})_n$ on $C^0([0,T]; H)$ is tight.
Using a similar argument, 
since $W^{\bar s, \ell}(0,T; V_1)$ is compactly embedded in
$C^0([0,T]; H)$, one can also show that the family of laws of
\[
  G_n(\varphi_{\lambda,n})\cdot W:=\int_0^\cdot G_n(\varphi_{\lambda,n}(s))\,\d W(s)
\]
is tight on $C^0([0,T]; H)$.
Moreover, taking into account the remarks in Subsection~\ref{ssec:not}
we identify $W$ with
a constant sequence of random variables with values in
$C^0([0,T]; U_1)$. Hence, we deduce that
the family of laws of
$(\varphi_{\lambda,n}, G_n(\varphi_{\lambda,n})\cdot W, W)_n$
is tight on the product space
\[
  C^0([0,T]; H) \times C^0([0,T]; H) \times C^0([0,T]; U_1)\,.
\]

By Prokhorov and Skorokhod theorems
(see \cite[Thm.~2.7]{ike-wata} and \cite[Thm.~1.10.4, Add.~1.10.5]{vaa-well})
and their weaker version by Jakubowski-Skorokhod (see e.~g.~\cite[Thm.~2.7.1]{breit-feir-hof}),
recalling the estimates \eqref{est1_n}--\eqref{est3_n}
there exists a probability space $(\tilde\Omega,\tilde\cF, \tilde\P)$ and measurable 
maps $\Lambda_n:(\tilde\Omega,\tilde\cF)\to(\Omega,\cF)$ such that 
$\tilde\P\circ \Lambda_n^{-1}=\P$ for every $n\in\enne$ and 
\begin{align*}
  \tilde\varphi_{\lambda,n}:=\varphi_{\lambda,n}\circ\Lambda_n\to \tilde\varphi_\lambda
  \qquad&\text{in } L^p(\tilde\Omega;C^0([0,T]; H)) \quad\,\forall\,p<\ell\,,\\
  \tilde\varphi_{\lambda,n}\wstarto \tilde\varphi_\lambda \qquad
  &\text{in } L^\ell_w(\tilde\Omega; L^\infty(0,T; V_1))
  \cap L^\ell(\tilde\Omega; W^{\bar s,\ell}(0,T; V_1^*))\,,\\
  \tilde\mu_{\lambda,n}:=\mu_{\lambda,n}\circ\Lambda_n \wto\tilde\mu_{\lambda}
  \qquad&\text{in } L^{\ell/2}(\tilde\Omega; L^2(0,T; V_1))\,,\\
  \nabla\tilde\mu_{\lambda,n}\wto\nabla\tilde\mu_{\lambda}
  \qquad&\text{in } L^{\ell}(\tilde\Omega; L^2(0,T; H))\,,\\
  \tilde I_{\lambda,n}:=(G_n(\varphi_{\lambda,n})\cdot W)\circ\Lambda_n \to \tilde I_{\lambda}
  \qquad&\text{in } L^p(\tilde\Omega; C^0([0,T]; H)) \quad\forall\,p<\ell\,,\\
  \tilde W_{n}:=W\circ\Lambda_n \to \tilde W 
  \qquad&\text{in } L^p(\tilde\Omega; C^0([0,T]; U_1)) \quad\forall\,p<\ell\,,
\end{align*}
for some measurable processes 
\begin{align*}
  &\tilde \varphi_\lambda \in L^\ell(\tilde\Omega; C^0([0,T]; H))\cap 
  L^\ell_w(\tilde\Omega; L^\infty(0,T; V_1))
  \cap L^\ell(\tilde\Omega; W^{\bar s,\ell}(0,T; V_1^*))\,,\\
  &\tilde\mu_\lambda \in L^{\ell/2}(\tilde\Omega; L^2(0,T; V_1))\,, 
  \quad\nabla\tilde\mu_\lambda\in L^\ell(\tilde\Omega; L^2(0,T; H))\,,\\
  &\tilde I_\lambda \in L^\ell(\tilde\Omega; C^0([0,T]; H))\,,\\
  & \tilde W \in L^\ell(\tilde \Omega; C^0([0,T]; U_1))\,. 
\end{align*}
Note that possibly enlarging the new probability space, it is not restrictive to 
suppose that $(\tilde\Omega,\tilde\cF,\tilde\P)$ is independent of $\lambda$.
Now, since $F_\lambda'$ is Lipschitz-continuous, we readily have
\[
F_\lambda'(\tilde\varphi_{\lambda,n})\to F_\lambda'(\tilde\varphi_\lambda) \qquad
\text{in } L^\ell(\tilde\Omega; L^2(0,T; H))\,,
\]
and similarly, since $G:H\to\cL^2(U,H)$ is Lipschitz-continuous,
\begin{align*}
  &\norm{G_n(\tilde\varphi_{\lambda,n})-
  G(\tilde\varphi_\lambda)}_{L^p(\tilde\Omega; L^2(0,T; \cL^2(U,H)))}\\
  &\leq\norm{G_n(\tilde\varphi_{\lambda,n})-
  G_n(\tilde\varphi_\lambda)}_{L^p(\tilde\Omega; L^2(0,T; \cL^2(U,H)))}
  +\norm{G_n(\tilde\varphi_\lambda)
  -G(\tilde\varphi_\lambda)}_{L^p(\tilde\Omega; L^2(0,T; \cL^2(U,H)))}\\
  &\leq\norm{G(\tilde\varphi_{\lambda,n})
  -G(\tilde\varphi_\lambda)}_{L^p(\tilde\Omega; L^2(0,T; \cL^2(U,H)))}
  +\norm{G_n(\tilde\varphi_\lambda)
  -G(\tilde\varphi_\lambda)}_{L^p(\tilde\Omega; L^2(0,T; \cL^2(U,H)))}\to 0\,,
\end{align*}
so that 
\[
  G_n(\tilde\varphi_{\lambda,n}) \to G(\tilde\varphi_\lambda) 
  \qquad\text{in } L^p(\tilde\Omega; L^2(0,T; \cL^2(U,H)))
  \quad\forall\,p<\ell\,.
\]
Moreover, since $\varphi_0\in V_1$ we also have that
$\varphi_0^n\to\varphi_0$ in $V_1$. 

Let us handle now the stochastic integral and identify the 
limit term $\tilde I_\lambda$:
we follow the classical arguments in \cite{fland-gat} and \cite[\S~8.4]{dapratozab}.
Introducing the filtration
\[
  \tilde\cF_{\lambda,n,t}:=\sigma\{\tilde\varphi_{\lambda,n}(s), 
  \tilde I_{\lambda,n}(s), \tilde W_n(s):s\in[0,t]\}\,, \qquad
  t\in[0,T]\,,
\]
clearly $\tilde W_n$ is adapted. Moreover,
since the maps $\Lambda_n$ preserve the laws,
we have that $\tilde W_n$ is a $Q_1$--Wiener process on $U_1$
(hence a cylindrical Wiener process on $U$), and 
$\tilde I_{\lambda,n}$ is the $H$-valued martingale 
\[
  \tilde I_{\lambda,n}(t)=
  \int_0^tG_n(\tilde \varphi_{\lambda,n}(s))\,\d \tilde W_n(s)\qquad\forall\,t\in[0,T]\,.
\]
These statements follow directly by the fact the $\Lambda_n$ preserves the laws,
and by definition of $Q_1$--Wiener process and stochastic integral
(see for example \cite[\S~4.5]{sap-witt-zimm} and \cite[\S~4.3]{vall-zimm}).

The next step is to show that the limit process $\tilde W$ is actually a 
cylindrical Wiener process in $U$, i.e.~that is a $Q_1$--Wiener process on $U_1$.
To this end, we introduce the filtration 
\[
  \tilde\cF_{\lambda,t}:=\sigma\{\tilde\varphi_{\lambda}(s),
  \tilde I_{\lambda}(s), \tilde W(s):s\in[0,t]\}\,, 
  \qquad
  t\in[0,T]\,,
\]
so that clearly $\tilde W$ is adapted to it. Also, it holds 
$\tilde W(0)=0$ since $\tilde W_n(0)=0$ and
$\tilde W_{n}\to \tilde W$ in $C^0([0,T]; U_1)$ as $n\to\infty$.
Moreover, for every $s,t\in[0,T]$ with $s\leq t$ and for 
every bounded continuous function 
$\psi\in C^0_b(C^0([0,s]; H)\times C^0([0,s]; H)\times C^0([0,s]; U_1))$, 
since $\tilde W_n$ is a $(\tilde\cF_{\lambda,n,t})_t$--martingale we have 
\[
  \tilde\E\left[(\tilde W_n(t) - \tilde W_n(s))
  \psi(\tilde \varphi_{\lambda,n}, \tilde I_{\lambda,n}, \tilde W_n)\right]=0
  \qquad\forall\,n\in\enne\,.
\]
Hence, letting $n\to\infty$ we deduce, thanks to the strong convergences above,
the continuity and boundedness of $\psi$, and the dominated convergence theorem, that 
\[
  \tilde\E\left[(\tilde W(t) - \tilde W(s))
  \psi(\tilde\varphi_{\lambda}, \tilde I_{\lambda}, \tilde W)\right]=0\,,
\]
yielding that $\tilde W$ is a $(\tilde\cF_{\lambda,t})_t$--martingale in $U_1$.
Similarly, since $\tilde W_n$ is a $Q_1$--Wiener process on $U_1$,
by the dominated convergence theorem it holds that, for every $h,k\in U_1$,
\begin{align*}
  0&=\tilde\E\left[\left((\tilde W_n(t), h)_{U_1}(\tilde W_n(t), k)_{U_1}
  - (\tilde W_n(s), h)_{U_1}(\tilde W_n(s), k)_{U_1}
  - (t-s)(Q_1h, k)_{U_1}\right)
  \psi(\tilde\varphi_{\lambda,n}, \tilde I_{\lambda,n}, \tilde W_n)\right]\\
  &\overset{n\to\infty}\longrightarrow
  \tilde\E\left[\left((\tilde W(t), h)_{U_1}(\tilde W(t), k)_{U_1}
  - (\tilde W(s), h)_{U_1}(\tilde W(s), k)_{U_1}
  - (t-s)(Q_1h, k)_{U_1}\right)
  \psi(\tilde\varphi_{\lambda}, \tilde I_{\lambda}, \tilde W)\right]\,,
\end{align*}
from which we get that the tensor quadratic variation of $\tilde W$ on $U_1$ is
\[
  \langle\!\langle\tilde W\rangle\!\rangle(t)=tQ_1\,, \qquad t\in[0,T]\,.
\]
Hence, by \cite[Thm.~4.6]{dapratozab} we deduce that $\tilde W$
is a $Q_1$--Wiener process on $U_1$, i.e.~a cylindrical Wiener process on $U$,
with respect to the filtration $(\tilde\cF_{\lambda,t})_t$.

Next, let us argue on the same same line for $\tilde I_\lambda$.
First of all, it is clear that $\tilde I_\lambda$ is $(\tilde \cF_{\lambda,t})_t$--adapted 
and $\tilde I_\lambda(0)=0$, since $\tilde I_{\lambda,n}(0)=0$ for all $n\in\enne$.
Secondly, for every $s,t\in[0,T]$ with $s\leq t$ and for every 
$\psi\in C^0_b(C^0([0,s]; H)\times C^0([0,s]; H)\times C^0([0,s]; U_1))$,
the martingale property of $\tilde I_{\lambda,n}$ and the dominated convergence theorem 
yield again
\begin{align*}
  0=
  \tilde\E\left[(\tilde I_{\lambda,n}(t) - \tilde I_{\lambda,n}(s))
  \psi(\tilde \varphi_{\lambda,n}, \tilde I_{\lambda,n}, \tilde W_n)\right]
  &\overset{n\to\infty}\longrightarrow
  \tilde\E\left[(\tilde I_{\lambda}(t) - \tilde I_{\lambda}(s))
  \psi(\tilde \varphi_{\lambda}, \tilde I_{\lambda}, \tilde W)\right]\,.
\end{align*}
We deduce that $\tilde I_\lambda$ is an $H$--valued $(\tilde\cF_{\lambda,t})_t$--martingale.
Similarly, since $\tilde I_{\lambda,n}=G_n(\tilde \varphi_{\lambda,n})\cdot \tilde W_n$,
we have the tensor quadratic variation 
\[
  \langle\!\langle\tilde I_{\lambda,n}\rangle\!\rangle(t)=\int_0^t G_n(\tilde \varphi_{\lambda,n}(s))\circ
  G_n(\tilde \varphi_{\lambda,n}(s))^*\,\d s\,, \qquad t\in[0,T]\,.
\]
Consequently, given arbitrary $h,k\in H$, 
the dominated convergence theorem yields again that 
\begin{align*}
  0&=
  \tilde\E\bigg[\Big((\tilde I_{\lambda,n}(t),h)_H(\tilde I_{\lambda,n}(t),k)_H 
  - (\tilde I_{\lambda,n}(s),h)_H(\tilde I_{\lambda,n}(s),k)_H \\
  &\qquad\qquad
  -\int_s^t(G_n(\tilde \varphi_{\lambda,n}(r))G_n(\tilde \varphi_{\lambda,n}(r))^*h, k)_H\,\d r
  \Big)
  \psi(\tilde\varphi_{\lambda,n}, \tilde I_{\lambda,n}, \tilde W_n)\bigg]\\
  &\overset{n\to\infty}\longrightarrow
  \tilde\E\bigg[\Big((\tilde I_{\lambda}(t),h)_H(\tilde I_{\lambda}(t),k)_H 
  - (\tilde I_{\lambda}(s),h)_H(\tilde I_{\lambda}(s),k)_H \\
  &\qquad\qquad
  -\int_s^t(G(\tilde \varphi_{\lambda}(r))G(\tilde \varphi_{\lambda}(r))^*h, k)_H\,\d r
  \Big)
  \psi(\tilde\varphi_{\lambda}, \tilde I_{\lambda}, \tilde W)\bigg]\,,
\end{align*}
so that the tensor quadratic variation of $\tilde I_\lambda$ is given by
\beq
  \label{cov1}
   \langle\!\langle\tilde I_{\lambda}\rangle\!\rangle(t)=
   \int_0^t G(\tilde \varphi_{\lambda}(s))\circ
  G(\tilde \varphi_{\lambda}(s))^*\,\d s\,, \qquad t\in[0,T]\,.
\eeq

Taking these remarks into account and setting 
$\tilde M_\lambda:= G(\tilde \varphi_\lambda)\cdot\tilde W$,
with the information collected so far we know that
$\tilde M_\lambda$ and $\tilde I_\lambda$ are $H$--valued martingales 
with respect to $(\tilde\cF_{\lambda,t})_t$ with tensor quadratic variations given by 
\beq
  \label{cov2}
  \langle\!\langle\tilde M_\lambda\rangle\!\rangle(t)
  =\langle\!\langle\tilde I_\lambda\rangle\!\rangle(t)=
  \int_0^t G(\tilde \varphi_{\lambda}(s))\circ
  G(\tilde \varphi_{\lambda}(s))^*\,\d s\,, \qquad t\in[0,T]\,.
\eeq
In order to conclude that actually $\tilde M_\lambda=\tilde I_\lambda$, we need
to exploit their tensor quadratic covariation. To this end, 
recalling the definition of $Q_1$ in Subsection~\ref{ssec:not},
by \cite[Thm.~3.2, p.~12]{Pard} we have that
\begin{align*}
  \langle\!\langle\tilde I_{\lambda,n}, \tilde W_n\rangle\!\rangle(t)&=
  \int_0^tG_n(\tilde\varphi_{\lambda,n}(s))\circ\iota^{-1}\,
  \d\langle\!\langle\tilde W_n\rangle\!\rangle(s)=
  \int_0^tG_n(\tilde \varphi_{\lambda,n}(s))\circ\iota^{-1}\circ Q_1\,\d s\\
  &=\int_0^tG_n(\tilde \varphi_{\lambda,n}(s))\circ\iota^*\,\d s\,, \qquad t\in[0,T]\,,
\end{align*}
yielding, taking adjoints (and recalling that $\iota:U\to U_1$ is the usual inclusion),
\[
  \langle\!\langle\tilde W_n, \tilde I_{\lambda,n}\rangle\!\rangle(t) = 
  \int_0^t\iota\circ G_n(\tilde\varphi_{\lambda,n}(s))^*\,\d s\,, \qquad t\in[0,T]\,.
\]
Hence, for every $s,t\in[0,T]$ with $s\leq t$,
for every $\psi\in C^0_b(C^0([0,s]; H)\times C^0([0,s]; H)\times C^0([0,s]; U_1))$,
for every $h\in U_1$ and $k\in H$,
we have again by the dominated convergence theorem 
\begin{align*}
  0&=\tilde\E\bigg[\Big((\tilde W_n(t), h)_{U_1}(\tilde I_{\lambda,n}(t), k)_H
  -(\tilde W_n(s), h)_{U_1}(\tilde I_{\lambda,n}(s), k)_H\\
  &\qquad\qquad-\int_s^t(G_n(\tilde\varphi_{\lambda,n}(r))^*k, h)_{U_1}\,\d r\Big)
  \psi(\tilde \varphi_{\lambda,n}, \tilde I_{\lambda, n}, \tilde W_n)\bigg]\\
  &\overset{n\to\infty}\longrightarrow
  \tilde\E\bigg[\Big((\tilde W(t), h)_{U_1}(\tilde I_{\lambda}(t), k)_H
  -(\tilde W(s), h)_{U_1}(\tilde I_{\lambda}(s), k)_H
  -\int_s^t(G(\tilde \varphi_{\lambda}(r))^*k, h)_{U_1}\,\d r\Big)
  \psi(\tilde \varphi_{\lambda}, \tilde I_{\lambda}, \tilde W)\bigg]\,.
\end{align*}
Consequently, we have that 
\beq
  \label{cov3}
  \langle\!\langle \tilde W, \tilde I_\lambda\rangle\!\rangle(t)=
  \int_0^t\iota\circ G(\tilde\varphi_{\lambda}(s))^*\,\d s\,, \qquad t\in[0,T]\,.
\eeq

We are now ready to conclude: indeed, taking \eqref{cov1}--\eqref{cov3} into account,
we get
\begin{align*}
  \langle\!\langle\tilde M_\lambda - \tilde I_\lambda\rangle\!\rangle&=
  \langle\!\langle\tilde M_\lambda\rangle\!\rangle + 
  \langle\!\langle\tilde I_\lambda\rangle\!\rangle - 
  2\langle\!\langle\tilde M_\lambda, \tilde I_\lambda\rangle\!\rangle\\
  &=2\int_0^\cdot G(\tilde\varphi_{\lambda}(s))\circ G(\tilde\varphi_{\lambda}(s))^*\,\d s
  - 2 \int_0^\cdot G(\tilde\varphi_{\lambda}(s))\,\d
  \langle\!\langle\tilde W, \tilde I_\lambda\rangle\!\rangle(s) = 0\,,
\end{align*}
which yields that 
\beq
  \label{id_stoc}
  \tilde I_\lambda=\tilde M_\lambda=\int_0^\cdot G(\tilde\varphi_\lambda(s))\,\d\tilde W(s)\,.
\eeq

Eventually, 
testing \eqref{eq1_app} by 
arbitrary $v\in V_1$ and integrating in time yields, after letting
$n\to\infty$ and taking \eqref{id_stoc} into account
\beq\label{eq1_app_lam}
\begin{split}
  \int_\OO\tilde\varphi_\lambda(t,x)v(x)\,\d x 
  &+ \int_{Q_t}m(\tilde\varphi_\lambda(s,x))
  \nabla\tilde\mu_\lambda(s,x)\cdot\nabla v(x)\,\d x\,\d s \\
  &=\int_\OO\varphi_0(x)v(x)\,\d x + 
  \int_\OO\left(\int_0^tG(\tilde\varphi_\lambda(s))\,\d \tilde W(s)\right)\!\!(x)v(x)\,\d x
\end{split}\eeq
for every $t\in[0,T]$, $\tilde\P$-almost surely. Indeed, this follows directly
form the convergences proved above, 
the fact that $m_n(\tilde\varphi_{\lambda,n})\to m(\tilde\varphi_\lambda)$
almost everywhere, the fact that
$|m_n|\leq m^*$,
and the dominated convergence theorem. Moreover,
testing \eqref{eq2_app} by $v\in V_1$ gives 
\beq\label{eq2_app_lam}
  \int_\OO\tilde\mu_\lambda(t,x)v(x)\,\d x = 
  \int_\OO\nabla\tilde\varphi_\lambda(s,x)\cdot\nabla v(x)\,\d x
  +\int_\OO F_\lambda'(\tilde\varphi_\lambda(t,x))v(x)\,\d x
\eeq
for almost every $t\in(0,T)$, $\tilde\P$-almost surely.

Finally, if
we take expectations
in \eqref{ito_app}, noting the stochastic integral 
is a martingale and that $\tilde\P\circ\Lambda_n^{-1}=\P$ for every $n$, 
we obtain
\begin{align*}
    \frac12\tilde\E\norm{\nabla\tilde\varphi_{\lambda,n}(t)}_H^2
    &+\tilde\E\norm{F_\lambda(\tilde\varphi_{\lambda,n}(t))}_{L^1(\OO)}
    +\tilde\E\int_{Q_t}m_n(\tilde\varphi_{\lambda,n}(s,x))
    |\nabla\tilde\mu_{\lambda,n}(s,x)|^2\,\d x\,\d s\\
    &\leq\frac12\norm{\nabla\varphi_0^{n}}_H^2 + \norm{F_\lambda(\varphi_0^{n})}_{L^1(\OO)}
    +\frac{C_G}{2}\tilde\E\int_0^t\norm{\nabla\tilde\varphi_{\lambda,n}(s)}_H^2\,\d s\\
    &+\frac12\tilde\E\int_0^t\sum_{k=0}^\infty\int_\OO
    F_\lambda''(\tilde\varphi_{\lambda,n}(s,x))|g_k(\tilde\varphi_{\lambda,n}(s,x))|^2\,\d x\,\d s
\end{align*}
for every $t\in[0,T]$. We want to let $n\to\infty$
using the convergences
proved above.
To this end, the first two terms on the left-hand side
and all the terms on the right-hand side
pass to the limit by weak lower semicontinuity and
the dominated convergence theorem
(recall that $F_\lambda''$ is continuous and bounded).
In order to pass to the limit 
by lower semicontinuity in the third term on the left-hand side,
it is sufficient to show that 
\beq\label{weak_conv}
\sqrt{m_n(\tilde\varphi_{\lambda,n})}\nabla\tilde\mu_{\lambda,n}\wto
\sqrt{m(\tilde\varphi_{\lambda})}\nabla\tilde\mu_{\lambda}
\quad\text{in } L^1(\tilde\Omega\times Q)\,.
\eeq
To prove this, note that
since $m_n(\tilde\varphi_{\lambda,n})\to m(\tilde\varphi_{\lambda})$
a.e.~in $\tilde\Omega\times Q$, for any arbitrary fixed $\sigma>0$, 
by the Severini-Egorov theorem there
is a measurable set $A_\sigma\subset \tilde\Omega\times Q$ such that 
$|A_\sigma^c|\leq\sigma$ and $m_n(\tilde\varphi_{\lambda,n})\to m(\tilde\varphi_\lambda)$
uniformly in $A_\sigma$. In particular, we have that 
$\sqrt{m_n(\tilde\varphi_{\lambda,n})}1_{A_\sigma}\to 
\sqrt{m(\tilde\varphi_\lambda)}1_{A_\sigma}$ in $L^\infty(\tilde\Omega\times Q)$. 
Consequently, for any $\zeta\in L^\infty(\tilde\Omega\times Q)^d$ we have
\[
  \int_{\tilde\Omega\times Q} \sqrt{m_n(\tilde\varphi_{\lambda,n})}
  \nabla\tilde\mu_{\lambda,n}\cdot\zeta=
  \int_{\tilde\Omega\times Q} 1_{A_\sigma}
  \sqrt{m_n(\tilde\varphi_{\lambda,n})}\nabla\tilde\mu_{\lambda,n}\cdot\zeta
  +\int_{A_\sigma^c} \sqrt{m_n(\tilde\varphi_{\lambda,n})}\nabla\tilde\mu_{\lambda,n}\cdot\zeta\,.
\]
Since $\sqrt{m_n(\tilde\varphi_{\lambda,n})}1_{A_\sigma}\to 
\sqrt{m(\tilde\varphi_\lambda)}1_{A_\sigma}$ in $L^\infty(\tilde\Omega\times Q)$ and 
$\zeta\in L^\infty(\tilde\Omega\times Q)$ we have
\[
  \int_{\tilde\Omega\times Q} 
  1_{A_\sigma}\sqrt{m_n(\tilde\varphi_{\lambda,n})}\nabla\tilde\mu_{\lambda,n}\cdot\zeta
  \to \int_{\tilde\Omega\times Q} 
  1_{A_\sigma}\sqrt{m(\tilde\varphi_{\lambda})}\nabla\tilde\mu_{\lambda}\cdot\zeta\,,
\]
while the H\"older inequality and the boundedness of $\nabla\tilde\mu_{\lambda,n}$ 
in $L^2(\tilde\Omega\times Q)$ yields
\[
\left|\int_{A_\sigma^c} 
\sqrt{m_n(\tilde\varphi_{\lambda,n})}\nabla\tilde\mu_{\lambda,n}\cdot\zeta\right|\leq
(m^*)^{1/2}
\norm{\nabla\tilde\mu_{\lambda,n}}_{L^2(\tilde\Omega\times Q)}
\norm{\zeta}_{L^2(A_\sigma^c)}
\leq c \norm{\zeta}_{L^\infty(\tilde\Omega\times Q)} \sigma^{1/2}\,,
\]
where $c>0$ is independent of $n$ and $\sigma$. Since $\sigma$ and $\zeta$ are arbitrary, 
we infer that \eqref{weak_conv} holds.
Hence, passing to the limit as $n\to\infty$ yields
by weak lower semicontinuity, for every $t\in[0,T]$,
\beq\label{en_ineq_aux}
  \begin{split}
    \frac12\tilde\E\norm{\nabla\tilde\varphi_{\lambda}(t)}_H^2
    &+\tilde\E\norm{F_\lambda(\tilde\varphi_{\lambda}(t))}_{L^1(\OO)}
    +\tilde\E\int_{Q_t}m(\tilde\varphi_{\lambda}(s,x))|\nabla\tilde\mu_{\lambda}(s,x)|^2\,\d x\,\d s\\
    &\leq\frac12\norm{\nabla\varphi_0}_H^2 + \norm{F_\lambda(\varphi_0)}_{L^1(\OO)}
    +\frac{C_G}{2}\tilde\E\int_0^t\norm{\nabla\tilde\varphi_{\lambda}(s)}_H^2\,\d s\\
    &+\frac12\tilde\E\int_0^t\sum_{k=0}^\infty\int_\OO
    F_\lambda''(\tilde\varphi_{\lambda}(s,x))|g_k(\tilde\varphi_{\lambda}(s,x))|^2\,\d x\,\d s\,.
    \end{split}
\eeq

\subsection{Uniform estimates in $\lambda$}
We prove here uniform estimates independently of $\lambda$.

First of all, since $\tilde\P\circ\Lambda_n^{-1}=\P$ for every $n\in\enne_+$,
from \eqref{est_mean} and weak lower semicontinuity it follows that
\[
  \norm{(\tilde\varphi_\lambda)_\OO}_{L^\ell(\tilde\Omega; C^0([0,T]))}\leq c
\]
for a certain $c>0$ independent of $\lambda$.

Secondly, from the estimate \eqref{ito_aux} we infer that 
\begin{align*}
   &\tilde\E\sup_{s\in[0,t]}\norm{\nabla\tilde\varphi_{\lambda,n}(s)}_H^{\ell} +
   \tilde\E\sup_{s\in[0,t]}\norm{F_\lambda(\tilde\varphi_{\lambda,n}(s))}_{L^1(\OO)}^{\ell/2}\\
   &\qquad+\tilde\E\sup_{s\in[0,t]}|(\tilde\mu_{\lambda,n}(s))_\OO|^{\ell/2}
   +\tilde\E\norm{\nabla\tilde\mu_{\lambda,n}}_{L^2(0,t; H)}^{\ell}\\
   &\leq c\left(1 + \tilde\E\norm{(\tilde\mu_{\lambda,n})_\OO}_{L^2(0,t)}^{\ell/2}
   + \tilde\E\norm{\nabla\tilde\varphi_{\lambda,n}}^\ell_{L^2(0,t; H)}\right)
   +\norm{F_\lambda(\varphi_0^{n})}_{L^1(\OO)}^{\ell/2}\\
   &\qquad+\frac12\tilde\E\left(\int_0^t\sum_{k=0}^\infty\int_\OO
  |F_\lambda''(\tilde\varphi_{\lambda,n}(s,x))|
  |G_{n}(\tilde\varphi_{\lambda,n}(s))u_k|^2(x)\,\d x\,\d s\right)^{\ell/2}\,,
\end{align*}
where again $c$ is independent of $\lambda$. We want to let $n\to\infty$
using weak lower semicontinuity of the norms at both sides.
To this end, since $\varphi_0^n\to\varphi_0$ in $V_1$
and $F_\lambda$ is bounded by a quadratic function by \eqref{Flam_quad},
$F_\lambda(\varphi_0^n)\to F_\lambda(\varphi_0)$ in $L^1(\OO)$ as $n\to\infty$.
Moreover, recalling also the strong convergences
$G_n(\tilde\varphi_{\lambda,n})\to G(\tilde\varphi_\lambda)$
in $L^p(\tilde\Omega; L^2(0,T; \cL^2(U,H)))$ and 
$\tilde\varphi_{\lambda,n}\to\tilde\varphi_\lambda$ in $L^p(\tilde\Omega; C^0([0,T]; H))$
for every $p<\ell$, proved in the previous subsection, we have in particular that
\[
  |F_\lambda''(\tilde\varphi_{\lambda,n})||G_{n}(\tilde\varphi_{\lambda,n})u_k|^2
  \to |F_\lambda''(\tilde\varphi_{\lambda})||G(\tilde\varphi_{\lambda})u_k|^2
  \quad\text{a.e.~in } \tilde\Omega\times(0,T)\times\OO\,,\quad\forall\,k\in\enne\,.
\]
Since $|F_\lambda''|\leq c_\lambda$, by {\bf ND3} we have
\[
  \int_\OO|F_\lambda''(\tilde\varphi_{\lambda,n})||G_{n}(\tilde\varphi_{\lambda,n})u_k|^2\leq
  c_\lambda\norm{G_n(\tilde\varphi_{\lambda,n})u_k}_H^2\leq 
  c_\lambda\norm{G(\tilde\varphi_{\lambda,n})u_k}_H^2\leq
  c_\lambda|\OO|\norm{g_k}^2_{L^\infty(\erre)}\,,
\]
and the dominated convergence theorem yields
\begin{align*}
  &\tilde\E\left(\int_0^t\sum_{k=0}^\infty\int_\OO
  |F_\lambda''(\tilde\varphi_{\lambda,n}(s,x))|
  |G_{n}(\tilde\varphi_{\lambda,n}(s))u_k|^2(x)\,\d x\,\d s\right)^{\ell/2}\\
  &\qquad\to
  \tilde\E\left(\int_0^t\sum_{k=0}^\infty\int_\OO
  |F_\lambda''(\tilde\varphi_{\lambda}(s,x))|
  |G(\tilde\varphi_{\lambda}(s))u_k|^2(x)\,\d x\,\d s\right)^{\ell/2}\,.
\end{align*}

We infer then, letting $n\to\infty$ and using weak lower semicontinuity, that
\begin{align*}
   &\tilde\E\sup_{s\in[0,t]}\norm{\nabla\tilde\varphi_{\lambda}(s)}_H^{\ell} +
   \tilde\E\sup_{s\in[0,t]}\norm{F_\lambda(\tilde\varphi_{\lambda}(s))}_{L^1(\OO)}^{\ell/2}\\
   &\qquad+\tilde\E\sup_{s\in[0,t]}|(\tilde\mu_{\lambda}(s))_\OO|^{\ell/2}
   +\tilde\E\norm{\nabla\mu_{\lambda}}_{L^2(0,t; H)}^{\ell}\\
   &\leq c\left(1 + \tilde\E\norm{(\tilde\mu_{\lambda})_\OO}_{L^2(0,t)}^{\ell/2}
   + \tilde\E\norm{\nabla\tilde\varphi_{\lambda}}^\ell_{L^2(0,t; H)}\right)
   +\norm{F_\lambda(\varphi_0)}_{L^1(\OO)}^{\ell/2}\\
   &\qquad+\frac12\tilde\E\left(\int_0^t\sum_{k=0}^\infty\int_\OO
  |F_\lambda''(\tilde\varphi_{\lambda}(s,x))|
  |G(\tilde\varphi_{\lambda}(s))u_k|^2(x)\,\d x\,\d s\right)^{\ell/2}\,,
\end{align*}
where the constant $c$ is independent of $\lambda$.
Let us bound the last two terms on the right-hand side uniformly in $\lambda$.
First of all, recalling that $F$ is a quadratic perturbation
of the convex function $\hat\gamma$, 
by definition of $F_\lambda$ we have that
$\norm{F_\lambda(\varphi_0)}_{L^1(\OO)}\leq
\norm{F(\varphi_0)}_{L^1(\OO)}$.
Secondly, using the H\"older inequality and assumption {\bf ND3} we have
\begin{align*}
  &\int_0^t\sum_{k=0}^\infty\int_\OO
  |F_\lambda''(\tilde\varphi_{\lambda}(s,x))|
  |G(\tilde\varphi_{\lambda}(s))u_k|^2(x)\,\d x\,\d s\\
  &\leq
  \int_0^t\norm{F_\lambda''(\tilde\varphi_\lambda(s))}_{L^1(\OO)}\,\d s
  \sum_{k=0}^\infty\norm{g_k(\tilde\varphi_\lambda)}_{L^\infty(Q)}^2\leq
  |Q|C_G\int_0^t\norm{F_\lambda''(\tilde\varphi_\lambda(s))}_{L^1(\OO)}\,\d s\,.
\end{align*}

Putting this information together, using the growth assumption on $F''$ in {\bf ND1}, we have then
\begin{align*}
   &\tilde\E\sup_{s\in[0,t]}\norm{\nabla\tilde\varphi_{\lambda}(s)}_H^{\ell} +
   \tilde\E\sup_{s\in[0,t]}\norm{F_\lambda(\tilde\varphi_{\lambda}(s))}_{L^1(\OO)}^{\ell/2}
   +\tilde\E\sup_{s\in[0,t]}|(\tilde\mu_{\lambda}(s))_\OO|^{\ell/2}
   +\tilde\E\norm{\nabla\tilde\mu_{\lambda}}_{L^2(0,t; H)}^{\ell}\\
   &\leq c\left(1 + \tilde\E\norm{(\tilde\mu_{\lambda})_\OO}_{L^2(0,t)}^{\ell/2}
   + \tilde\E\norm{\nabla\tilde\varphi_{\lambda}}^\ell_{L^2(0,t; H)}
   +\tilde\E\norm{F_\lambda(\tilde\varphi_\lambda)}_{L^1(0,t; L^1(\OO))}^{\ell/2}\right)
   +\norm{F(\varphi_0)}_{L^1(\OO)}^{\ell/2}
\end{align*}
where the constant $c$ (possibly updated) is independent of $\lambda$.
The Gronwall lemma yields then
the estimates
\begin{align}
  \label{est1_lam}
  \norm{\tilde\varphi_{\lambda}}_{L^\ell(\tilde\Omega; L^\infty(0,T; V_1))}&\leq c\,,\\
  \label{est2_lam}
  \norm{\tilde\mu_{\lambda}}_{L^{\ell/2}(\tilde\Omega; L^2(0,T; V_1))} + 
  \norm{\nabla\tilde\mu_{\lambda}}_{L^\ell(\tilde\Omega; L^2(0,T; H))} &\leq c\,.
\end{align}
Thanks to {\bf ND3} we deduce that 
\[
  \norm{G(\tilde\varphi_{\lambda})}_{L^\infty(\tilde\Omega\times(0,T); \cL^2(U,H))
  \cap L^\ell(\tilde\Omega; L^\infty(0,T; \cL^2(U,V_1)))}\leq c\,,
\]
hence also by \cite[Lem.~2.1]{fland-gat}, for any $s\in(0,1/2)$,
\beq
\label{est3_lam}
  \norm{\int_0^\cdot G(\tilde\varphi_{\lambda}(s))
  \,\d \tilde W(s)}_{L^\ell(\tilde\Omega; W^{s,\ell}(0,T; V_1))}\leq c_{s}\,.
\eeq
By comparison in \eqref{eq1_app_lam} we infer that 
\beq\label{est4_lam}
  \norm{\tilde \varphi_{\lambda}}_{L^\ell(\tilde\Omega; W^{\bar s,\ell}(0,T; V_1^*))}\leq c\,,
\eeq
where again $c$ is independent of $\lambda$, and 
$\bar s\in(1/\ell,1/2)$ is fixed. Finally, testing the variational equation \eqref{eq2_app_lam}
by $\gamma_\lambda(\tilde\varphi_\lambda)=
F_\lambda'(\tilde\varphi_\lambda)+C_F \tilde\varphi_\lambda$
and rearranging the terms
we have, for almost every $t\in(0,T)$,
\begin{align*}
  &\int_\OO\gamma_\lambda'(\tilde\varphi_\lambda(t,x))|\nabla\tilde\varphi_\lambda(t,x)|^2\,\d x
  +\int_\OO|F_\lambda'(\tilde\varphi_\lambda(t,x))|^2\,\d x\\
  &=\int_\OO
  \tilde\mu_\lambda(t,x)F_\lambda'(\tilde\varphi_\lambda(t,x))\,\d x 
  +C_F\int_\OO \tilde\mu_\lambda(t,x)\tilde\varphi_\lambda(t,x)\,\d x 
  -C_F\int_\OO F_\lambda'(\tilde\varphi_\lambda(t,x))\tilde\varphi_\lambda(t,x)\,\d x\,.
\end{align*}
Since the first term on the left-hand side is nonnegative by monotonicity of $\gamma_\lambda$,
the Young inequality yields, after integrating in time,
\[
  \norm{F_\lambda'(\tilde\varphi_\lambda)}_{L^2(0,T; H)}^2\leq
  \frac12\norm{F_\lambda'(\tilde\varphi_\lambda)}_{L^2(0,T; H)}^2
  +\frac32\norm{\tilde\mu_\lambda}_{L^2(0,T; H)}^2 
  +\frac32 C_F^2\norm{\tilde\varphi_\lambda}_{L^2(0,T; H)}^2\,,
\]
so that by \eqref{est1_lam}--\eqref{est2_lam} we have
\beq
  \label{est5_lam}
  \norm{F_\lambda'(\tilde\varphi_\lambda)}_{L^{\ell/2}(\tilde\Omega; L^2(0,T; H))}\leq c\,.
\eeq
By comparison in \eqref{eq2_app_lam} we deduce that 
$\Delta\tilde\varphi_\lambda\in L^{\ell/2}(\tilde\Omega; L^2(0,T; H))$ and by elliptic regularity
\beq\label{est6_lam}
  \norm{\tilde\varphi_\lambda}_{L^{\ell/2}(\tilde\Omega; L^2(0,T; V_2))}\leq c\,.
\eeq

\subsection{Passage to the limit as $\lambda\to0$}
\label{ssec:lim_lambda}
We perform here the passage to the limit as $\lambda\to 0$
and recover a martingale solution to the original problem \eqref{eq1}--\eqref{eq4}.
Since the arguments are very similar to the ones in Subsection~\ref{ssec:lim_n},
we shall omit the details.

Since we have the compact inclusion 
$L^\infty(0,T; V_1)\cap W^{\bar s,\ell}(0,T; V_1^*)\cembed C^0([0,T]; H)$,
using the estimates \eqref{est1_lam} and \eqref{est4_lam}
and arguing exactly as in Subsection~\ref{ssec:lim_n},
one readily deduce that the family of laws of $(\tilde\varphi_\lambda)_\lambda$
on $C^0([0,T]; H)$ is tight. Similarly, estimate \eqref{est3_lam}
and the compact inclusion $W^{\bar s,\ell}(0,T; V_1)\cembed C^0([0,T]; H)$
yields that the family of laws of $(G(\tilde\varphi_\lambda)\cdot \tilde W)_\lambda$
is tight on $C^0([0,T]; H)$.
In particular, the family of laws of
$(\tilde\varphi_{\lambda}, G(\tilde\varphi_{\lambda})\cdot \tilde W, \tilde W)_\lambda$
is tight on the product space
\[
  C^0([0,T]; H) \times C^0([0,T]; H) \times C^0([0,T]; U_1)\,.
\]

By Prokhorov and Jakubowski-Skorokhod theorems
(see again the references 
\cite[Thm.~2.7]{ike-wata}, \cite[Thm.~1.10.4, Add.~1.10.5]{vaa-well}, and
\cite[Thm.~2.7.1]{breit-feir-hof}),
there exists a further probability space $(\hat\Omega,\hat\cF,\hat\P)$
and measurable 
maps $\Xi_\lambda:(\hat\Omega,\hat\cF)\to(\tilde\Omega,\tilde\cF)$ such that 
$\hat\P\circ \Xi_\lambda^{-1}=\tilde\P$ for every $\lambda>0$ and 
\begin{align*}
  \hat\varphi_{\lambda}:=\tilde\varphi_{\lambda}\circ\Xi_\lambda\to \hat\varphi
  \qquad&\text{in } L^p(\hat\Omega;C^0([0,T]; H)) \quad\,\forall\,p<\ell\,,\\
  \hat\varphi_{\lambda}\wstarto \hat\varphi 
  \qquad&\text{in } L^\ell_w(\hat\Omega; L^\infty(0,T; V_1))\,,\\
  \hat\varphi_{\lambda}\wto \hat\varphi \qquad&\text{in } L^{\ell/2}(\hat\Omega; L^2(0,T; V_2))
  \cap L^\ell(\hat\Omega; W^{\bar s,\ell}(0,T; V_1^*))\,,\\
  \hat\mu_{\lambda}:=\tilde\mu_{\lambda}\circ\Xi_\lambda \wto\hat\mu
  \qquad&\text{in } L^{\ell/2}(\hat\Omega; L^2(0,T; V_1))\,,\\
  \nabla\hat\mu_{\lambda} \wto\nabla\hat\mu
  \qquad&\text{in } L^{\ell}(\hat\Omega; L^2(0,T; H))\,,\\
  F_\lambda'(\hat\varphi_\lambda)\wto \hat\xi
  \qquad&\text{in } L^{\ell/2}(\hat\Omega; L^2(0,T; H))\,,\\
  \hat I_{\lambda}:=(G(\tilde\varphi_{\lambda})\cdot \tilde W)\circ\Xi_\lambda \to \hat I
  \qquad&\text{in } L^p(\hat\Omega; C^0([0,T]; H)) \quad\forall\,p<\ell\,,\\
  \hat W_\lambda:=\tilde W\circ\Xi_\lambda \to \hat W 
  \qquad&\text{in } L^p(\hat \Omega; C^0([0,T]; U_1)) \quad\forall\,p<\ell\,,
\end{align*}
for some measurable processes 
\begin{align*}
  &\hat\varphi\in L^\ell(\hat\Omega; C^0([0,T]; H))\cap L^\ell_w(\hat \Omega; L^\infty(0,T; V_1))
  \cap L^{\ell/2}(\hat\Omega; L^2(0,T; V_2))
  \cap L^\ell(\hat\Omega; W^{\bar s,\ell}(0,T; V_1^*))\,,\\
  &\hat\mu \in L^{\ell/2}(\hat\Omega; L^2(0,T; V_1))\,, 
  \quad\nabla\hat\mu\in L^\ell(\hat\Omega; L^2(0,T; H))\,,\\
  &\hat\xi\in  L^{\ell/2}(\hat\Omega; L^2(0,T; H))\,,\\
  &\hat I \in L^\ell(\hat \Omega; C^0([0,T]; H))\,,\\
  &\hat W \in L^\ell(\hat\Omega; C^0([0,T]; U_1))\,.
\end{align*}
Since $G:H\to\cL^2(U,H)$ is Lipschitz-continuous, we also have 
\[
  G(\hat\varphi_{\lambda}) \to G(\hat\varphi) \qquad\text{in } L^p(\hat\Omega; L^2(0,T; \cL^2(U,H)))
  \quad\forall\,p<\ell\,.
\]
Moreover, since $F'=\gamma - C_F I$ and $\gamma$ is maximal monotone,
by the strong-weak closure of maximal montone operators (see \cite[Lem.~2.3]{barbu-monot})
and the strong convergence of $(\hat\varphi_\lambda)_\lambda$,
we have that 
\[
  \xi = F'(\hat\varphi)  \quad\text{a.e.~in } \hat\Omega\times(0,T)\times\OO\,.
\]

Following the exact same argument of Subsection~\ref{ssec:lim_n}, 
we introduce the filtration
\[
  \hat\cF_{\lambda,t}:=
  \sigma\{\hat\varphi_{\lambda}(s), \hat I_{\lambda}(s), \hat W_\lambda(s):s\in[0,t]\}\,, 
  \qquad t\in[0,T]\,,
\]
and infer that $\hat I_{\lambda}$ is the $H$-valued martingale  given by
\[
  \hat I_{\lambda}(t)=\int_0^tG(\hat\varphi_{\lambda}(s))\,
  \d \hat W_\lambda(s)\qquad\forall\,t\in[0,T]\,.
\]
Now, arguing again as in Subsection~\ref{ssec:lim_n}, thanks to
the strong convergences of $\hat\varphi_{\lambda}\to\hat\varphi$ and 
$G(\hat\varphi_{\lambda})\to G(\hat\varphi)$ proved above,
we can suitably enlarge the
probability space
$(\hat\Omega, \hat\cF, \hat\P)$ and find a
saturated and right-continuous filtration
$(\hat\cF_t)_{t\in[0,T]}$ such that 
$\hat W$ is a 
cylindrical Wiener process on the stochastic basis
$(\hat\Omega, \hat\cF, (\hat\cF_t)_{t\in[0,T]}, \hat\P)$
and
\[
  \hat I(t)=\int_0^tG(\hat\varphi(s))\,\d \hat W(s) \qquad \forall\,t\in[0,T]\,.
\]

Now, since $\hat\P\circ\Xi_\lambda=\tilde \P$ for every $\lambda>0$,
from \eqref{eq1_app_lam} and \eqref{eq2_app_lam} it follows that
\begin{align*}
  \int_\OO\hat\varphi_\lambda(t,x)v(x)\,\d x 
  &+ \int_{Q_t}m(\hat\varphi_\lambda(s,x))\nabla\hat\mu_\lambda(s,x)\cdot\nabla v(x)\,\d x\,\d s \\
  &=\int_\OO\varphi_0(x)v(x)\,\d x + 
  \int_\OO\left(\int_0^tG(\hat\varphi_\lambda(s))\,\d \hat W_\lambda(s)\right)\!\!(x)v(x)\,\d x
\end{align*}
for every $t\in[0,T]$, $\hat\P$-almost surely, and
\[
  \int_\OO\hat\mu_\lambda(t,x)v(x)\,\d x = 
  \int_\OO\nabla\hat\varphi_\lambda(s,x)\cdot\nabla v(x)\,\d x
  +\int_\OO F_\lambda'(\hat\varphi_\lambda(t,x))v(x)\,\d x
\]
for almost every $t\in(0,T)$, $\hat\P$-almost surely. Hence, 
using the convergences proved above, 
the continuity and boundedness of $m$ together with the 
dominated convergence theorem, we 
can let $\lambda\to 0$ in the variational formulations and
obtain exactly \eqref{var-for-ND}, and $\hat\mu=-\Delta\hat\varphi+F'(\hat\varphi)$.

In order to prove the energy inequality \eqref{en_ineq}
we note that since $\hat\P\circ\Xi_\lambda^{-1}=\tilde\P$
for all $\lambda>0$, from
\eqref{en_ineq_aux} we infer that 
\beq\label{en_ineq_aux2}\begin{split}
    \frac12\hat\E\norm{\nabla\hat\varphi_{\lambda}(t)}_H^2
    &+\hat\E\norm{F_\lambda(\hat\varphi_{\lambda}(t))}_{L^1(\OO)}
    +\hat\E\int_{Q_t}m(\hat\varphi_{\lambda}(s,x))|\nabla\hat\mu_{\lambda}(s,x)|^2\,\d x\,\d s\\
    &\leq\frac12\norm{\nabla\varphi_0}_H^2 
    + \norm{F_\lambda(\varphi_0)}_{L^1(\OO)}
    +\frac{C_G}{2}\hat\E\int_0^t\norm{\nabla\hat\varphi_{\lambda}(s)}_H^2\,\d s\\
    &+\frac12\hat\E\int_0^t\sum_{k=0}^\infty\int_\OO
    F_\lambda''(\hat\varphi_{\lambda}(s,x))|g_k(\hat\varphi_{\lambda}(s,x))|^2\,\d x\,\d s
\end{split}\eeq
for every $t\in[0,T]$. 
We want to let again $\lambda\to0$ and use
the convergences just proved.
To this end, for the second term on the left-hand side
note that $F_\lambda$ is a 
quadratic perturbation of the convex function $\hat\gamma_\lambda$, where
$\hat\gamma_\lambda(\hat\varphi_\lambda)\geq\hat\gamma(J_\lambda\hat\varphi_\lambda)$
and $J_\lambda:=(I+\lambda\gamma)^{-1}:\erre\to\erre$ is the resolvent of $\gamma$.
Noting that 
\[
  |J_\lambda\hat\varphi_\lambda - \hat\varphi|\leq |J_\lambda\hat\varphi_\lambda - \hat\varphi_\lambda|
  +|\hat\varphi_\lambda - \hat\varphi|\leq \lambda |\gamma_\lambda(\hat\varphi_\lambda)|
  +|\hat\varphi_\lambda - \hat\varphi|\,,
\]
we easily infer that $J_\lambda\hat\varphi_\lambda\to\hat\varphi$ a.e.~in $\hat\Omega\times Q$,
so that by weak lower semicontinuity and the Fatou lemma
\[
  \hat\E\norm{F(\hat\varphi)}_{L^1(\OO)}\leq\liminf_{\lambda\to0}
  \hat\E\norm{F_\lambda(\hat\varphi_\lambda)}_{L^1(\OO)}\,.
\]
Moreover, for the third term on the left-hand side of \eqref{en_ineq_aux2},
arguing exactly as in the proof of
\eqref{weak_conv} we have that 
\[
  \hat\E\int_{Q_t}m(\hat\varphi(s,x))|\nabla\hat\mu(s,x)|^2\,\d x\,\d s\leq
  \liminf_{\lambda\to 0}
  \hat\E\int_{Q_t}m(\hat\varphi_\lambda(s,x))|\nabla\hat\mu_\lambda(s,x)|^2\,\d x\,\d s\,.
\]
Furthermore, for the second term on the right-hand side of \eqref{en_ineq_aux2}
note that $F_\lambda(\varphi_0)\leq F(\varphi_0)$. Finally,
the last term on the right-hand side of \eqref{en_ineq_aux2} can pass to the limit
by the dominated convergence theorem.
Indeed, since
$F_\lambda' = \gamma_\lambda - C_FI=\gamma\circ J_\lambda - C_FI$, we have
$F_\lambda''=(\gamma'\circ J_\lambda)J_\lambda' - C_F$, from which 
\[
  |F_\lambda''(\hat\varphi_{\lambda})||g_k(\hat\varphi_{\lambda})|^2\leq
  \norm{g_k}^2_{L^\infty(\erre)}|F_\lambda''(\hat\varphi_{\lambda})|=
  \norm{g_k}^2_{L^\infty(\erre)}
  |\gamma'(J_\lambda(\hat\varphi_\lambda))J_\lambda'(\hat\varphi_\lambda)-C_F|\,.
\]
Recalling that $J_\lambda:\erre\to\erre$ is $1$-Lipschitz continuous, 
it holds $|J_\lambda'|\leq1$ and we deduce that
\[
  |F_\lambda''(\hat\varphi_{\lambda})||g_k(\hat\varphi_{\lambda})|^2
  \leq \norm{g_k}^2_{L^\infty(\erre)}\left(|\gamma'(J_\lambda(\hat\varphi_\lambda))| + C_F\right)\,,
\]
where, by definition of $\gamma$  we have 
$\gamma'(J_\lambda(\hat\varphi_\lambda))=
F''(J_\lambda\hat\varphi_\lambda) + C_F$.
From these computations and assumption {\bf ND1} it follows that
\begin{align*}
  |F_\lambda''(\hat\varphi_{\lambda})||g_k(\hat\varphi_{\lambda})|^2
  &\leq \norm{g_k}^2_{L^\infty(\erre)}\left(|F''(J_\lambda(\hat\varphi_\lambda))|
  + 2C_F\right)
  \leq c \norm{g_k}^2_{L^\infty(\erre)} \left(1 + |F(J_\lambda(\hat\varphi_\lambda))|\right)\,.
\end{align*}
The term in brackets on the right-hand side is uniformly 
integrable in $\hat\Omega\times Q$ because
$F(J_\lambda(\hat\varphi_\lambda))\to F(\hat\varphi)$ in
$L^1(\hat\Omega\times Q)$,
hence so is the left-hand side by comparison. 
Hence, recalling that 
$\sum_{k=0}^\infty\norm{g_k}_{L^\infty(\erre)}^2\leq C_G$, 
the Vitali convergence theorem yields
\[
    \hat\E\int_0^t\sum_{k=0}^\infty\int_\OO
    F_\lambda''(\hat\varphi_{\lambda}(s,x))|g_k(\hat\varphi_{\lambda}(s,x))|^2\,\d x\,\d s
    \to
    \hat\E\int_0^t\sum_{k=0}^\infty\int_\OO
    F''(\hat\varphi(s,x))|g_k(\hat\varphi(s,x))|^2\,\d x\,\d s\,.
\]
Letting then $\lambda\to0$ in \eqref{en_ineq_aux2}
taking into account these remarks yields exactly,
by weak lower semicontinuity, the energy inequality \eqref{en_ineq}.

In order to conclude we only need to prove the last assertion of 
Theorem~\ref{thm1}. To this end, note that if
$d\geq3$ and $|F''(r)|\leq C_F(1+|r|^{\frac{2}{d-2}})$
for all $r\in\erre$,
recalling that $V_1\embed L^{\frac{2d}{d-2}}(\OO)$
and noting that $\frac{2-d}{2d} + \frac{1}{d} = \frac12$,
by H\"older inequality we have
\begin{align*}
  \norm{\nabla F_\lambda'(\hat\varphi_\lambda)}_H&=
  \norm{F_\lambda''(\hat\varphi_\lambda)\nabla\hat\varphi_\lambda}_H\leq
  C_F\norm{(1+|\hat\varphi_\lambda|^{\frac{2}{d-2}})\nabla\hat\varphi_\lambda}_H\\
  &\leq C_F\norm{1 + |\hat\varphi_\lambda|^{\frac{2}{d-2}}}_{L^d(\OO)}
  \norm{\nabla\hat\varphi_\lambda}_{L^{\frac{2d}{d-2}}(\OO)}
  \leq C_F\norm{1 + \hat\varphi_\lambda}^{\frac{2}{d-2}}_{L^{\frac{2d}{d-2}}(\OO)}
  \norm{\nabla\hat\varphi_\lambda}_{L^{\frac{2d}{d-2}}(\OO)}\\
  &\leq c\norm{1+\hat\varphi_\lambda}^{\frac{2}{d-2}}_{V_1}\norm{\hat\varphi_\lambda}_{V_2}\,,
\end{align*}
from which it follows that 
\[
  \norm{\nabla F_\lambda'(\hat\varphi_\lambda)}_{L^2(0,T; H)}\leq
  c\norm{1+\hat\varphi_\lambda}^{\frac{2}{d-2}}_{L^\infty(0,T; V_1)}\norm{\hat\varphi_\lambda}_{L^2(0,T; V_2)}\,.
\]
Using \eqref{est1_lam} and \eqref{est6_lam}, since $\frac{2}{d-2}\leq2$ the right-hand side
is uniformly bounded in $L^{\ell/4}(\hat\Omega)$, so that by weak lower
semicontinuity we have $F'(\hat\varphi)\in L^{\ell/4}(\hat\Omega; L^2(0,T; V_1))$.
Since $\hat\mu=-\Delta\hat\varphi + F'(\hat\varphi)$, we conclude by 
elliptic regularity. 
If $d=2$, the same argument works
using the embedding $V_1\embed L^q(\OO)$ for
all $q\in\mathopen[2,+\infty\mathclose)$.
This concludes the proof of Theorem~\ref{thm1}.

%%%%%%%%%%%%%%%%%%%%%%%%%%%%%%%%%%%%%%%%%%%%%%%%%%%%%

\section{Degenerate mobility and irregular potential}
\label{sec:D}

This section is devoted to proving Theorem~\ref{thm2}.
The main idea of the proof is the following.
We approximate the irregular potential $F$ and the mobility
$m$ using a suitable regularization, depending
on a parameter $\eps>0$, introduced
in \cite{ell-gar} in the deterministic setting.
We show that the $\eps$-approximated problem 
admits martingale solutions thanks to the already
proved Theorem~\ref{thm2}.
Finally, exploiting the compatibility assumptions
between $F$, $m$, and $G$, we prove uniform
estimates on the solutions and pass to the limit by 
monotonicity and stochastic compactness arguments.

Let us also mention that a similar approximation in $\eps>0$
was used in \cite{FGGS2019}, in the study of a nonlocal 
Cahn--Hilliard--Navier--Stokes deterministic model with degenerate mobility.
Here the authors prove that the system with $\eps>0$ fixed
admits a solution by relying on a time--discretisation argument, and then
they show convergence as $\eps\to0$.
In our case, the idea is similar, but existence of solution at $\eps$ fixed
is obtained using the $\lambda$--approximation of the {\bf ND} case
instead of the time--discretisation.
This choice was meant to avoid technical time--measurability issues in the stochastic setting;
nonetheless, we point out that
both techniques are feasible also in the stochastic case.

\subsection{The approximation}
For brevity of notation, we shall denote by 
$\eps$ a fixed real sequence $(\eps_n)_n\subset(0,1/4)$
converging to $0$ as $n\to\infty$.
We shall briefly say that $\eps\to0$.

Since $F_2\in C^2([-1,1])$
we can extend
it to the whole $\erre$ as
\[
  \tilde F_2:\erre\to\erre\,,\qquad
  \tilde F_2(r):=
  \begin{cases}
  F_2(r) \quad&\text{if } |r|\leq1\,,\\
  F_2(-1) + F_2'(-1)(r+1) + \frac12F''(-1)(r+1)^2 \quad&\text{if } r<-1\,,\\
  F_2(1) + F_2'(1)(r-1) + \frac12F''(1)(r-1)^2 \quad&\text{if } r>1\,,
  \end{cases}
\]
so that $\tilde F_{2}\in C^2(\erre)$ and $\|\tilde F_{2}''\|_{C^0(\erre)}\leq\norm{F_2''}_{C^0([-1,1])}$.
Moreover, 
for every $\eps\in(0,1/4)$ we define
$F_{1,\eps}:\erre\to\erre$ as the unique 
function of class $C^2$ such that 
$F_{1,\eps}(0)=F_1(0)$, $F_{1,\eps}'(0)=F_1'(0)$, and
\[
  F''_{1,\eps}(r)=
  \begin{cases}
  F_1''(r) \quad&\text{if } |r|\leq 1-\eps\,,\\
  F_1''(-1+\eps) \quad&\text{if } r<-1+\eps\,,\\
  F_1''(1-\eps) \quad&\text{if } r>1-\eps\,.
  \end{cases}
\]
With this notation, we introduce
the regularized potential $F_\eps:=F_{1,\eps} + \tilde F_2$
and note that, by {\bf D1}, we have 
$F_\eps \in C^2(\erre)$ with $F_\eps''\in L^\infty(\erre)$.
In particular, $F_\eps$ satisfies {\bf ND1} and \eqref{extra_growth}.
Moreover, by definition we have that $F_\eps=F$ on $[-1+\eps, 1-\eps]$.

We define the approximated mobility
\[
  m_\eps: \erre\to\erre\,, \qquad
  m_\eps(r):=
  \begin{cases}
  m(r) \quad&\text{if } |r|\leq 1-\eps\,,\\
  m(-1+\eps) \quad&\text{if } r<-1+\eps\,,\\
  m(1-\eps) \quad&\text{if } r>1-\eps\,.
  \end{cases}
\]
Note that by assumption {\bf D2} we have that 
$m_\eps\in C^0(\erre)$ with
\[
0<\min\{m(-1+\eps), m(1-\eps)\}\leq m_\eps(r)\leq \norm{m}_{C^0([-1,1])} \quad\forall\,r\in\erre\,,
\]
so that $m_\eps$ satisfies {\bf ND2} and $m_\eps=m$ on $[-1+\eps, 1-\eps]$.
Furthermore, we define $M_\eps$ as the unique
function in $C^2(\erre)$ such that 
\[
  M_\eps:\erre\to\mathopen[0,+\infty\mathclose)\,, \qquad
  M_\eps(0)=M_\eps'(0)=0\,, \qquad
  M''_\eps(r):=\frac1{m_\eps(r)}\,, \quad r\in\erre\,.
\]
In particular, note that by definition of $m_\eps$ and {\bf D2} we have $M_\eps''\in L^\infty(\erre)$.

We define the approximated operator $G_\eps:H\to\cL^2(U,H)$ setting
\[
  G_\eps(v)u_k:=g_{k,\eps}(v) \qquad v\in H\,, \quad k\in\enne\,,
\]
where, for every $k\in\enne$,
\[
  g_{k,\eps}:\erre\to\erre\,, \qquad
  g_{k,\eps}(r):=
  \begin{cases}
    g_k(r) \quad&\text{if } |r|\leq 1-\eps\,,\\
    g_k(-1+\eps) \quad&\text{if } r<-1+\eps\,,\\
    g_k(1-\eps) \quad&\text{if } r>1-\eps\,.
  \end{cases}
\]
Note that $g_{k,\eps}\in W^{1,\infty}(\erre)$ for every $k\in\enne$
and that 
\[
  C_{G_\eps}=\sum_{k=0}^\infty\norm{g_{k,\eps}}_{W^{1,\infty}(\erre)}^2\leq
  \sum_{k=0}^\infty\norm{g_{k}}_{W^{1,\infty}(-1,1)}^2 \leq L_G<+\infty\,,
\]
so that $G_\eps$ satisfies assumption {\bf ND3} and $G_\eps=G$ on $\mathcal B_{1-\eps}$.

We deduce that the assumptions {\bf ND1}--{\bf ND4}
are satisfied by the set of data 
$(F_\eps, m_\eps, G_\eps, \varphi_0)$:
hence, Theorem~\ref{thm1}
ensures that for every $\eps\in(0,1/4)$
there exists a martingale solution
\[
\left(\Omega_\eps, \cF_\eps, (\cF_{\eps,t})_{t\in[0,T]}, \P_\eps, W_\eps,\varphi_\eps, \mu_\eps\right)
\]
to the problem \eqref{eq1}--\eqref{eq4}
in the sense of Definition~\ref{def:sol-ND}.
Since $\eps\in(0,1/4)$ is chosen in
a countable set of points
(see the remark at the beginning of the subsection),
by suitably enlarging the probability 
spaces we shall suppose that 
such martingale solutions are of the form
\[
\left(\bar\Omega, \bar\cF, (\bar\cF_{t})_{t\in[0,T]}, \bar\P, \bar W, \varphi_\eps, \mu_\eps\right)\,,
\]
for a certain stochastic basis $(\bar\Omega, \bar\cF, (\bar\cF_{t})_{t\in[0,T]}, \bar\P)$.
Theorem~\ref{thm1} ensures also the regularities
\begin{align*}
  \varphi_\eps&\in L^\ell(\bar\Omega; C^0([0,T]; H))\cap L^\ell_w(\bar\Omega; L^\infty(0,T; V_1))\\
  &\quad\cap L^{\ell/2}(\bar\Omega; L^2(0,T; V_2))\cap 
  L^{\ell/4}(\bar\Omega; L^2(0,T; H^3(\OO)))\,,\\
  \mu_\eps &\in L^{\ell/2}(\bar\Omega; L^2(0,T; V_1))\,, 
  \quad \nabla\mu_\eps\in L^\ell(\bar\Omega; L^2(0,T; H^d))\,,\\
  F_\eps'(\varphi_\eps) &\in L^{\ell/2}(\bar\Omega; L^2(0,T; H))\cap 
  L^{\ell/4}(\bar\Omega; L^2(0,T; V_1))\,.
\end{align*}
for every $\ell\in\mathopen[2,+\infty\mathclose)$, and
the variational formulation \eqref{var-for-ND} reads
\beq\label{var_eps}\begin{split}
  \int_\OO\varphi_\eps(t,x)v(x)\,\d x 
  &+ \int_{Q_t}m_\eps(\varphi_\eps(s,x))\nabla\mu_\eps(s,x)\cdot\nabla v(x)\,\d x\,\d s \\
  &=\int_\OO\varphi_0(x)v(x)\,\d x + 
  \int_\OO\left(\int_0^tG_\eps(\hat\varphi_\eps(s))\,\d W'(s)\right)\!\!(x)v(x)\,\d x
\end{split}\eeq
for every $v\in V_1$,
for every $t\in[0,T]$, $\bar\P$-almost surely, where
$\mu_\eps=-\Delta\varphi_\eps + F'_\eps(\varphi_\eps)$.

\subsection{Uniform estimates in $\eps$}
We show here uniform estimates on the approximated solutions.

\noindent {\sc First estimate.} First of all, taking $v=1$ in \eqref{var_eps},
using It\^o's formula and the
Burkholder-Davis-Gundy and Young inequalities as
in the proof of \eqref{est_mean} yields 
\beq
  \label{est_mean_eps}
  \norm{(\varphi_\eps)_\OO}_{L^\ell(\bar\Omega; C^0([0,T]))}\leq c\,,
\eeq
for a positive constant $c$ independent of $\eps$.

Secondly, the energy inequality \eqref{en_ineq} implies that 
\begin{align*}
    \frac12\sup_{r\in[0,t]}\bar\E\norm{\nabla\varphi_\eps(r)}_H^2
    &+\sup_{r\in[0,t]}\bar\E\norm{F_\eps(\varphi_\eps(r))}_{L^1(\OO)}
    +\bar\E\int_{Q_t}m_\eps(\varphi_\eps(s,x))|\nabla\mu_\eps(s,x)|^2\,\d x\,\d s\\
    &\leq\frac12\norm{\nabla\varphi_0}_H^2 + \norm{F(\varphi_0)}_{L^1(\OO)}
    +\frac{C_{G_\eps}}{2}\bar\E\int_0^t\norm{\nabla\varphi_\eps(s)}_H^2\,\d s\\
    &+\frac12\bar\E\int_0^t\sum_{k=0}^\infty\int_\OO
    F''_\eps(\varphi_\eps(s,x))|g_{k_\eps}(\varphi_\eps(s,x))|^2\,\d x\,\d s\,.
\end{align*}
Now, we have already shown that $C_{G_\eps}\leq L_G$ for every $\eps$.
Moreover, by assumption {\bf D3} and the definitions of $F_\eps$ and $m_\eps$, we have that
for every $r\in[-1+\eps, 1-\eps]$
\[
  F_\eps''(r)|g_{k,\eps}(r)|^2=
  F''(r)|g_k(r)|^2\leq \norm{\sqrt{F''}g_k}^2_{L^\infty(-1,1)}\,,
\]
while for every $r<-1+\eps$
\begin{align*}
  F_\eps''(r)|g_k(r)|^2&=(F_{1}''(-1+\eps) + \tilde F''_2(r))|g_k(-1+\eps)|^2\\
  &=F''(-1+\eps)|g_k(-1+\eps)|^2 + (\tilde F''_2(r) - F_2''(-1+\eps))|g_k(-1+\eps)|^2\\
  &\leq\norm{\sqrt{F''}g_k}^2_{L^\infty(-1,1)} + 2\norm{F_2''}_{C^0([-1,1])}\norm{g_k}_{L^\infty(-1,1)}^2\,,
\end{align*}
and similarly for every $r>1-\eps$
\begin{align*}
  F_\eps''(r)|g_k(r)|^2&=(F_{1}''(1-\eps) + \tilde F''_2(r))|g_k(1-\eps)|^2\\
  &=F''(1-\eps)|g_k(1-\eps)|^2 + (\tilde F''_2(r) - F_2''(1-\eps))|g_k(1-\eps)|^2\\
  &\leq\norm{\sqrt{F''}g_k}^2_{L^\infty(-1,1)} + 2\norm{F_2''}_{C^0([-1,1])}\norm{g_k}_{L^\infty(-1,1)}^2\,.
\end{align*}
Hence, we deduce that 
\[
  F_\eps''(r)|g_{k,\eps}(r)|^2\leq
   \norm{\sqrt{F''}g_k}^2_{L^\infty(-1,1)} + 2\norm{F_2''}_{C^0([-1,1])}\norm{g_k}_{L^\infty(-1,1)}^2
   \qquad\forall\,r\in\erre\,.
\]
Substituting in the energy inequality and recalling the definition of $L_G$ yields then
\begin{align*}
    &\frac12\sup_{r\in[0,t]}\bar\E\norm{\nabla\varphi_\eps(r)}_H^2
    +\sup_{r\in[0,t]}\bar\E\norm{F_\eps(\varphi_\eps(r))}_{L^1(\OO)}
    +\bar\E\int_{Q_t}m_\eps(\varphi_\eps(s,x))|\nabla\mu_\eps(s,x)|^2\,\d x\,\d s\\
    &\leq\frac12\norm{\nabla\varphi_0}_H^2 + \norm{F(\varphi_0)}_{L^1(\OO)}
    +\frac{L_G}{2}\bar\E\int_0^t\norm{\nabla\varphi_\eps(s)}_H^2\,\d s
    +(1+ 2\norm{F_2}_{C^0([-1,1])}) L_G \frac{|Q|}{2}\,.
\end{align*}
The Gronwall lemma and estimate \eqref{est_mean_eps} imply then
\begin{align}
  \label{est1_eps}
  \norm{\varphi_\eps}_{C^0([0,T]; L^2(\bar\Omega; V_1))}&\leq c\,,\\
  \label{est2_eps}
  \norm{F_\eps(\varphi_\eps)}_{L^\infty(0,T; L^1(\bar\Omega\times\OO))}&\leq c\,,\\
  \label{est3_eps}
  \norm{\sqrt{m_\eps(\varphi_\eps)}\nabla\mu_\eps}_{L^2(\bar\Omega; L^2(0,T; H))}&\leq c\,.
\end{align}
Noting that $\norm{g_{k,\eps}}_{W^{1,\infty}(\erre)}\leq\norm{g_k}_{W^{1,\infty}(-1,1)}$
for every $k\in\enne$, 
thanks to assumption {\bf D3} and \eqref{est1_eps} we also have that 
\[
  \norm{G_\eps(\varphi_\eps)}_{L^\infty(\bar\Omega\times(0,T); \cL^2(U,H))\cap 
  L^\infty(0,T; L^2(\bar\Omega; \cL^2(U,V_1)))}\leq c\,,
\]
which implies by \cite[Lem.~2.1]{fland-gat}
that,
for every $s\in(0,1/2)$ and for every $p\in\mathopen[2,+\infty\mathclose)$,
\beq\label{est4_eps}
  \norm{\int_0^\cdot G_\eps(\varphi_\eps(s))\,\d W'(s)}_{L^p(\bar\Omega; W^{s,p}(0,T; H))
  \cap L^2(\bar\Omega; W^{s,2}(0,T; V_1))}\leq c_{s,p}\,.
\eeq
By comparison in \eqref{var_eps} we deduce from \eqref{est3_eps} and \eqref{est4_eps} that,
for every $s\in(0,1/2)$,
\beq
  \label{est5_eps}
  \norm{\varphi_\eps}_{L^2(\bar\Omega; W^{s,2}(0,T; V_1^*))}\leq c_s\,.
\eeq

\noindent {\sc Second estimate.}
The idea is now to write It\^o's formula for $\int_\OO M_\eps(\varphi_\eps)$. In order to do this,
note that $M_\eps\in C^2(\erre)$
and
$M_\eps''\in L^\infty(\erre)$, so that $M'_\eps$ is Lipschitz-continuous. 
In particular, we have that 
$v\mapsto M_\eps'(v)$ is well defined and continuous from $V_1$ to $V_1$.
Hence, we can apply It\^o's formula in the variational setting
\cite[Thm.~4.2]{Pard} and obtain
\begin{align*}
  \int_\OO M_\eps(\varphi_\eps(t,x))\,\d x &+ 
  \int_{Q_t}m_\eps(\varphi_\eps(s,x))M_\eps''(\varphi_\eps(s,x))
  \nabla\mu_\eps(s,x)\cdot\nabla\varphi_\eps(s,x)\,\d x\,\d s\\
  &=\int_\OO M_\eps(\varphi_0(x))\,\d x
  +\int_0^t\left(M_\eps'(\varphi_\eps(s)), G_\eps(\varphi_\eps(s))\,\d W'(s)\right)_H\\
  &+\frac12\int_0^t\sum_{k=0}^\infty\int_\OO
  M_\eps''(\varphi_\eps(s,x))|g_k(\varphi_\eps(s,x))|^2\,\d x\,\d s\,.
\end{align*}
Noting that $M_\eps'' m_\eps=1$ by definition of $M_\eps$
and recalling that 
$\mu_\eps=-\Delta\varphi_\eps + F_\eps'(\varphi_\eps)$,
since the stochastic integral on the right-hand side is a martingale
taking expectations yields
\begin{align*}
  &\bar\E\int_\OO M_\eps(\varphi_\eps(t,x))\,\d x +
  \bar\E\int_{Q_t}|\Delta\varphi_\eps(s,x)|^2\,\d x\,\d s
  +\bar\E\int_{Q_t}F_\eps''(\varphi_\eps(s,x))|\nabla\varphi_\eps(s,x)|^2\,\d x\,\d s\\
  &=\int_\OO M_\eps(\varphi_0(x))\,\d x + 
  \frac12\bar\E\int_0^t\sum_{k=0}^\infty\int_\OO
  M_\eps''(\varphi_\eps(s,x))|g_k(\varphi_\eps(s,x))|^2\,\d x\,\d s\,.
\end{align*}
Now, by definition of $M_\eps$ and $m_\eps$ we have $M_\eps(\varphi_0)\leq M(\varphi_0)$
almost everywhere in $\OO$. Moreover,
by assumption {\bf D3} and the definitions of $M_\eps$ and $G_\eps$, we have that
for every $r\in[-1+\eps, 1-\eps]$
\[
  M_\eps''(r)|g_{k,\eps}(r)|^2=
  M''(r)|g_k(r)|^2\leq \norm{\sqrt{M''}g_k}^2_{L^\infty(-1,1)}\,,
\]
while for every $r<-1+\eps$
\[
  M_\eps''(r)|g_k(r)|^2=
  M''(-1+\eps)|g_k(-1+\eps)|^2
  \leq\norm{\sqrt{M''}g_k}^2_{L^\infty(-1,1)}\,,
\]
and similarly for every $r>1-\eps$
\[
  M_\eps''(r)|g_k(r)|^2=
  M''(1-\eps)|g_k(1-\eps)|^2
  \leq\norm{\sqrt{M''}g_k}^2_{L^\infty(-1,1)}\,.
\]
Hence, we deduce that 
\[
  M_\eps''(r)|g_{k,\eps}(r)|^2\leq
   \norm{\sqrt{M''}g_k}^2_{L^\infty(-1,1)}
   \quad\forall\,r\in\erre\,.
\]
Taking into account these remarks and recalling that $F_\eps'=F_{1,\eps}'+\tilde F_2'$, we have then
\begin{align*}
  &\bar\E\int_\OO M_\eps(\varphi_\eps(t,x))\,\d x +
  \bar\E\int_{Q_t}|\Delta\varphi_\eps(s,x)|^2\,\d x\,\d s
  +\bar\E\int_{Q_t}F_{1,\eps}''(\varphi_\eps(s,x))|\nabla\varphi_\eps(s,x)|^2\,\d x\,\d s\\
  &\leq\int_\OO M(\varphi_0(x))\,\d x + 
  L_G\frac{|Q|}2 - \bar\E\int_{Q_t}\tilde F_2''(\varphi_\eps(s,x))|\nabla\varphi_\eps(s,x)|^2\,\d x\,\d s\,.
\end{align*}
Since $F_{1,\eps}''\geq0$ by definition and $\|\tilde F_2''\|_{L^\infty(\erre)}\leq\norm{F_2''}_{C^0([-1,1])}$,
recalling {\bf D4} and the estimate \eqref{est1_eps}
we deduce by elliptic regularity that
\begin{align}
  \label{est6_eps}
  \norm{M_\eps(\varphi_\eps)}_{C^0([0,T]; L^1(\bar\Omega\times\OO))}&\leq c\,,\\
  \label{est7_eps}
  \norm{\varphi_\eps}_{L^2(\bar\Omega; L^2(0,T; V_2))}&\leq c\,.
\end{align}

\noindent{\sc Third estimate.}
We prove now an estimate allowing to 
obtain some $L^\infty$-bounds on the limiting solution:
we are inspired here by some computations performed in \cite[p.~414]{ell-gar}
(see also \cite[\S~4.1.1]{frig-lam-roc}).
Note that by definition of $M_\eps$ and $m_\eps$ we have that,
for every $r>1$,
\begin{align*}
  M_\eps(r)&=M_\eps(1-\eps) + M_\eps'(1-\eps)(r-1+\eps) + \frac12M_\eps''(1-\eps)(r-1+\eps)^2\\
  &\geq\frac12M_\eps''(1-\eps)(r-1+\eps)^2 = \frac{(r-1+\eps)^2}{2m_\eps(1-\eps)}
  \geq\frac{(r-1)^2}{2m_\eps(1-\eps)}\,,
\end{align*}
and similarly, for every $r<-1$,
\begin{align*}
  M_\eps(r)&=M_\eps(-1+\eps) + M_\eps'(-1+\eps)(r+1-\eps) + \frac12M_\eps''(-1+\eps)(r+1-\eps)^2\\
  &\geq\frac12M_\eps''(-1+\eps)(r+1-\eps)^2 = \frac{(r+1-\eps)^2}{2m_\eps(-1+\eps)}
  \geq\frac{(r+1)^2}{2m_\eps(-1+\eps)}=\frac{(|r|-1)^2}{2m_\eps(-1+\eps)}\,.
\end{align*}
We infer that
\[
  (|r|-1)^2_+\leq 2M_\eps(r)\max\{m_\eps(1-\eps), m_\eps(-1+\eps)\} \quad\forall\,r\in\erre\,.
\]
Now, since $m_\eps(1-\eps)=m(1-\eps)$ and $m(1)=0$ we have
\[
  |m_\eps(1-\eps)|=|m(1-\eps)| = |m(1-\eps)-m(1)| \leq \norm{m'}_{L^\infty(-1,1)}\eps
\]
and similarly, since $m_\eps(-1+\eps)=m(-1+\eps)$ and $m(-1)=0$,
\[
  |m_\eps(-1+\eps)|=|m(-1+\eps)| = |m(-1+\eps)-m(-1)| \leq \norm{m'}_{L^\infty(-1,1)}\eps\,.
\]
We deduce then that
\[
  (|r|-1)^2_+\leq 2\eps \norm{m'}_{L^\infty(-1,1)}
  M_\eps(r) \quad\forall\,r\in\erre\,,
\]
yielding, together with \eqref{est6_eps},
\beq
  \label{est8_eps}
  \norm{(|\varphi_\eps|-1)_+}_{C^0([0,T];L^2(\bar\Omega; H))} \leq c\sqrt{\eps}\,.
\eeq

\subsection{Passage to the limit as $\eps\to0$}
The argument are similar to the ones performed in
Subsections~\ref{ssec:lim_n} and \ref{ssec:lim_lambda},
so we will omit the technical details for brevity. 

First of all, by \cite[Cor.~5, p.~86]{simon} we have the compact inclusion
\[
  L^2(0,T; V_2)\cap W^{s,2}(0,T; V_1^*)\cembed L^2(0,T; V_1)\,.
\]
Hence, using the estimates \eqref{est5_eps} and \eqref{est7_eps},
together with the Markov inequality, arguing as in 
Subsection~\ref{ssec:lim_n} we easily infer that the
laws of $(\varphi_\eps)_\eps$ on $L^2(0,T; V_1)$
are tight. Secondly, 
fixing $\bar s\in(0,1/2)$ and $\bar p\geq2$ such that
$\bar s\bar p>1$,
again by \cite[Cor.~5, p.~86]{simon} 
we have also the compact inclusion
\[
  W^{\bar s,\bar p}(0,T; H)\cap W^{\bar s,2}(0,T; V_1)\cembed C^0([0,T]; V_1^*)\cap L^2(0,T; H)\,,
\]
so that estimate \eqref{est4_eps} yields by the same argument that 
the sequence of laws of $(G_\eps(\varphi_\eps)\cdot \bar W)_\eps$
on $C^0([0,T]; V_1^*)\cap L^2(0,T; H)$ are tight.
In particular, the family of laws of
$(\varphi_{\eps}, G_\eps(\varphi_\eps)\cdot \bar W, \bar W)_\eps$
is tight on the product space
\[
  L^2(0,T; V_1)\times \left(C^0([0,T]; V_1^*)\cap L^2(0,T; H)\right) \times C^0([0,T]; U_1)\,.
\]

Recalling that $L^2(0,T; H)$, $L^2(0,T; V_1)$, and $L^2(0,T; V_2)$
endowed with their weak topologies are sub-Polish,
by Prokhorov and Jakubowski-Skorokhod theorems
(see again \cite[Thm.~2.7]{ike-wata}, \cite[Thm.~1.10.4, Add.~1.10.5]{vaa-well}, and
\cite[Thm.~2.7.1]{breit-feir-hof})
and the estimates \eqref{est1_eps}--\eqref{est7_eps},
there exists a 
probability space $(\hat\Omega,\hat\cF,\hat\P)$
and measurable 
maps $\Theta_\eps:(\hat\Omega,\hat\cF)\to(\bar\Omega,\bar\cF)$ such that 
$\hat\P\circ \Theta_\eps^{-1}=\bar\P$ for every $\eps\in(0,1/4)$ and 
\begin{align*}
  \hat\varphi_{\eps}:=\varphi_\eps\circ\Theta_\eps\to \hat\varphi
  \qquad&\text{in } L^2(0,T; V_1) \quad\hat\P\text{-a.s.}\,,\\
  \hat\varphi_{\eps}\wstarto \hat\varphi \qquad&\text{in } L^\infty(0,T; L^2(\hat\Omega; V_1))\,,\\
  \hat\varphi_{\eps}\wto \hat\varphi \qquad&\text{in } L^2(0,T; V_2)
  \cap W^{\bar s,2}(0,T; V_1^*)\quad\hat\P\text{-a.s.}\,,\\
  \hat\eta_{\eps}:=(m_\eps(\varphi_\eps)\nabla\mu_{\eps})\circ\Theta_\eps \wto\hat\eta
  \qquad&\text{in } L^2(0,T; H) \quad\hat\P\text{-a.s.}\,,\\
  \hat I_{\eps}:=(G_\eps(\varphi_{\eps})\cdot \bar W)\circ\Theta_\eps \to \hat I
  \qquad&\text{in } L^p(\hat\Omega; C^0([0,T]; V_1^*)\cap L^2(0,T; H)) 
  \quad\forall\,p\in\mathopen[2,+\infty\mathclose)\,,\\
  \hat W_\eps:=\bar W\circ\Theta_\eps \to \hat W 
  \qquad&\text{in } L^p(\hat \Omega; C^0([0,T]; U)) \quad\forall\,p\in\mathopen[2,+\infty\mathclose)\,,
\end{align*}
for some measurable processes 
\begin{align*}
  &\hat\varphi\in L^\infty(0,T; L^2(\hat\Omega; V_1))
  \cap L^{2}(\hat\Omega; L^2(0,T; V_2))
  \cap L^2(\hat\Omega; W^{\bar s,2}(0,T; V_1^*))\,,\\
  &\hat\eta \in L^2(\hat\Omega; L^2(0,T; H^d))\,,\\
  &\hat I \in L^p(\hat \Omega; C^0([0,T]; V^*)\cap L^2(0,T; H))
  \quad\forall\,p\in\mathopen[2,+\infty\mathclose)\,,\\
  &\hat W \in L^p(\hat\Omega; C^0([0,T]; U))
  \quad\forall\,p\in\mathopen[2,+\infty\mathclose)\,.
\end{align*}
Furthermore, by lower semicontinuity and estimate \eqref{est8_eps}
it follows that 
\[
  |\hat\varphi(t)|\leq 1 \quad\text{a.e.~in }\hat\Omega\times\OO\,,\quad\forall\,t\in[0,T]\,,
\]
while the estimates \eqref{est1_eps} yields also the convergence
\[
  \hat\varphi_{\eps}\to \hat\varphi
  \qquad\text{in } L^p(\hat\Omega; L^2(0,T; V_1)) \quad\forall\,p\in\mathopen[1,2\mathclose)\,.
\]
Now, by definition of $G_\eps$ and recalling that 
$\|g_{k,\eps}'\|_{L^\infty(\erre)}\leq
\norm{g_k'}_{L^\infty(-1,1)}$, we have that 
\begin{align*}
  \norm{G_\eps(\hat\varphi_\eps)-G(\hat\varphi)}_{\cL^2(U,H)}^2&\leq
  2\norm{G_\eps(\hat\varphi_\eps)-G_\eps(\hat\varphi)}_{\cL^2(U,H)}^2
  +2\norm{G_\eps(\hat\varphi)-G(\hat\varphi)}_{\cL^2(U,H)}^2\\
  &=2\sum_{k=0}^\infty\left(
  \norm{g_{k,\eps}(\hat\varphi_\eps)-g_{k,\eps}(\hat\varphi)}_H^2+
  \norm{g_{k,\eps}(\hat\varphi)-g_k(\hat\varphi)}_H^2\right)\\
  &\leq2\sum_{k=0}^\infty\left(
  \norm{g_{k,\eps}'}_{L^\infty(\erre)}^2\norm{\hat\varphi_\eps-\hat\varphi}_H^2
  +\norm{g_{k,\eps}(\hat\varphi)-g_k(\hat\varphi)}_H^2\right)\\
  &\leq 2L_G\norm{\hat\varphi_\eps-\hat\varphi}_H^2
  +2\sum_{k=0}^\infty\norm{g_{k,\eps}(\hat\varphi)-g_k(\hat\varphi)}_H^2\,,
\end{align*}
where
\begin{align*}
  &\sum_{k=0}^\infty\norm{g_{k,\eps}(\hat\varphi)-g_k(\hat\varphi)}_H^2\\
  &=\sum_{k=0}^\infty\left(\norm{(g_{k}(-1+\eps)-g_k(\hat\varphi))1_{\{\hat\varphi<-1+\eps\}}}_H^2
  +\norm{(g_{k}(1-\eps)-g_k(\hat\varphi))1_{\{\hat\varphi>1-\eps\}}}_H^2\right)\\
  &\leq\sum_{k=0}^\infty\norm{g_k'}^2_{L^\infty(-1,1)}\left(
  \norm{(\hat\varphi+1-\eps)1_{\{\hat\varphi<-1+\eps\}}}_H^2
  +\norm{(\hat\varphi-1+\eps)1_{\{\hat\varphi>1-\eps\}}}_H^2
  \right)\\
  &\leq L_G\left(
  \norm{(\hat\varphi+1-\eps)1_{\{\hat\varphi<-1+\eps\}}}_H^2
  +\norm{(\hat\varphi-1+\eps)1_{\{\hat\varphi>1-\eps\}}}_H^2
  \right)\,.
\end{align*}
Since $|\hat\varphi|\leq1$ almost everywhere,
we have that $1_{\{\hat\varphi<-1+\eps\}}\to0$ and 
$1_{\{\hat\varphi>1-\eps\}}\to 0$ almost everywhere, so that 
the dominated convergence theorem yields
\[
  G_\eps(\hat\varphi_\eps) \to G(\hat\varphi) \quad\text{in }
  L^p(\hat\Omega; L^2(0,T; \cL^2(U,H))) \quad\forall\,
  p\in\mathopen[1,2\mathclose)\,.
\]

Following now the same argument of Subsection~\ref{ssec:lim_n},
we define
the filtration
\[
  \hat\cF_{\eps,t}:=\sigma\{\hat\varphi_{\eps}(s), \hat I_{\eps}(s), \hat W_\eps(s):s\in[0,t]\}\,, 
  \qquad t\in[0,T]\,,
\]
and infer that
\[
  \hat I_{\eps}(t)=\int_0^tG_\eps(\hat\varphi_{\eps}(s))\,\d \hat W_\eps(s)\qquad\forall\,t\in[0,T]\,.
\]
Following again Subsection~\ref{ssec:lim_n}, thanks to
the strong convergences of $\hat\varphi_{\eps}\to\hat\varphi$ and 
$G_\eps(\hat\varphi_{\eps})\to G(\hat\varphi)$ proved above,
we deduce that, possibly enlarging the
probability space
$(\hat\Omega, \hat\cF, \hat\P)$, there is a
saturated and right-continuous filtration
$(\hat\cF_t)_{t\in[0,T]}$ such that 
$\hat W$ is a 
cylindrical Wiener process on the stochastic basis
$(\hat\Omega, \hat\cF, (\hat\cF_t)_{t\in[0,T]}, \hat\P)$
and
\[
  \hat I(t)=\int_0^tG(\hat\varphi(s))\,\d \hat W(s) \qquad \forall\,t\in[0,T]\,.
\]

Passing to the weak limit as $\eps\to0$ yields then, for every $v\in V_1$,
\beq\label{var_eta}\begin{split}
  \int_\OO\hat\varphi(t,x)v(x)\,\d x 
  &+ \int_{Q_t}\hat\eta(s,x)\cdot\nabla v(x)\,\d x\,\d s \\
  &=\int_\OO\varphi_0(x)v(x)\,\d x + 
  \int_\OO\left(\int_0^tG(\hat\varphi(s))\,\d \hat W(s)\right)\!\!(x)v(x)\,\d x
\end{split}\eeq
$\hat\P$-almost surely, for every $t\in[0,T]$.
It only remains to identify the limit $\hat\eta$.

Let us show that, for every $\zeta\in L^\infty(0,T; V_1\cap L^d(\OO))$, it holds
\beq\label{def_eta}\begin{split}
  \int_{Q} \hat\eta(s,x)\cdot\zeta(s,x)\,\d x\,\d s
  &=\int_{Q}\Delta\hat\varphi(s,x)\div\left[m_\eps(\hat\varphi(s,x))\zeta(s,x)\right]\,\d x\,\d s\\
  &+\int_{Q}m(\hat\varphi(s,x))F''(\hat\varphi(s,x))\nabla\hat\varphi(s,x)\cdot\zeta(s,x)\,\d x\,\d s
  \qquad\hat\P\text{-a.s.}\,.
\end{split}\eeq
First of all, note that all the terms in \eqref{def_eta} are well defined.
Indeed, by boundedness of $m_\eps$ we have
$m_\eps(\hat\varphi_\eps)\zeta\in L^\infty(0,T; H^d)$,
and recalling that $V_1\embed L^{\frac{2d}{d-2}}(\OO)$ and
$\frac{d-2}{2d} + \frac1d = \frac12$,
by H\"older inequality
\[
\div\left[m_\eps(\hat\varphi_\eps)\zeta\right] =
m_\eps'(\hat\varphi_\eps)\nabla\hat\varphi_\eps\cdot\zeta
+m_\eps(\hat\varphi_\eps)\div\zeta\,,
\]
with
\begin{align*}
  \norm{\div\left[m_\eps(\hat\varphi_\eps)\zeta\right]}_H&\leq 
  \norm{m'_\eps}_{L^{\infty}(\erre)}
  \norm{\nabla\varphi}_{L^{\frac{2d}{d-2}}}\norm{\zeta}_{L^d(\OO)}
  +\norm{m_\eps}_{L^{\infty}(\erre)}\norm{\div\zeta}_{H}\\
  &\leq c \norm{m}_{W^{1,\infty}(-1,1)}
  \left(\norm{\hat\varphi_\eps}_{V_2}\norm{\zeta}_{L^d(\OO)}
  +\norm{\eta}_{V_1}\right)\,.
\end{align*}
Hence $\div\left[m_\eps(\hat\varphi_\eps)\zeta\right] \in L^\infty(0,T; H)$.
Furthermore, the second term on the right-hand side is also
well defined since $m_\eps F_\eps''\in L^\infty(\erre)$.

By density, it is not restrictive to prove \eqref{def_eta}
for every $\zeta\in L^\infty(0,T;V_2\cap H^{2\bar m}(\OO))$,
where $\bar m$ is such that $H^{2\bar m}(\OO)\embed L^\infty(\OO)$.
Since $\hat\mu_\eps=-\Delta\hat\varphi_\eps + F_\eps'(\hat\varphi_\eps)$
and we have the regularity 
$\hat\varphi_\eps\in L^2(0,T; V_2\cap H^3(\OO))$,
integration by parts yields
\beq\label{def_eta_ep}\begin{split}
  &\int_{Q} m_\eps(\hat\varphi_\eps(s,x))\nabla\hat\mu_\eps(s,x)\cdot\zeta(s,x)\,\d x\,\d s\\
  &=-\int_{Q}\nabla\Delta\hat\varphi_\eps(s,x)\cdot m_\eps(\hat\varphi_\eps(s,x))\zeta(s,x)\,\d x\,\d s\\
  &\qquad+\int_Q
  F_\eps''(\hat\varphi_\eps(s,x))m_\eps(\hat\varphi_\eps(s,x))\nabla\hat\varphi_\eps(s,x)\cdot\zeta(s,x)\,\d x\,\d s\\
  &=\int_{Q} \Delta\hat\varphi_\eps(s,x)\div\left[m_\eps(\hat\varphi_\eps(s,x))\zeta(s,x)\right]\,\d x\,\d s\\
  &\qquad+\int_Q
  F_\eps''(\hat\varphi_\eps(s,x))m_\eps(\hat\varphi_\eps(s,x))\nabla\hat\varphi_\eps(s,x)\cdot\zeta(s,x)
  \,\d x\,\d s
\end{split}\eeq
so that we need to show that we can pass to the limit on the right-hand side.
Let us start from the first-term. This is of the form
\[
  \Delta\hat\varphi_\eps m_\eps'(\hat\varphi_\eps)\nabla\hat\varphi_\eps\cdot\zeta
  +\Delta\hat\varphi_\eps m_\eps(\hat\varphi_\eps)\div\zeta\,,
\]
where 
\beq\label{conv_aux1}
  \Delta\hat\varphi_\eps\wto\Delta\hat\varphi \quad\text{in } L^2(0,T; H)
  \quad\hat\P\text{-a.s.}
\eeq
Since $m_\eps(\hat\varphi_\eps)\to m(\hat\varphi)$
almost everywhere and $|m_\eps|\leq |m|$, 
the dominated convergence theorem
yields that 
\[
m_\eps(\hat\varphi_\eps)\to m(\hat\varphi) \quad\text{in }
L^{p}(Q) \quad\hat\P\text{-a.s.}\,,
\qquad\forall\,p\in\mathopen[2,+\infty\mathclose)\,.
\]
Since $\zeta\in L^\infty(0,T; H^{2\bar m}(\OO))$ and 
$H^{2\bar m}(\OO)\embed L^\infty(\OO)$, this implies in particular
\beq
  \label{conv_aux2}
  m_\eps(\hat\varphi_\eps)\div\zeta\to m(\hat\varphi)\div\zeta
  \quad\text{in } L^2(0,T; H) \quad\hat\P\text{-a.s.}\,.
\eeq
Moreover, since $m_\eps'=0$ on $(-\infty,-1+\eps)\cup(1-\eps,+\infty)$,
we have that
\[
m_\eps'(\hat\varphi_\eps)=m'(\hat\varphi_\eps)1_{\{|\hat\varphi_\eps|\leq 1-\eps\}}
\to m'(\hat\varphi) \quad\text{a.e.~in } \hat\Omega\times Q\,,
\]
hence also, thanks to the strong convergence
$\hat\varphi_\eps\to\hat\varphi$ in $L^2(0,T; V_1)$,
that
\[
m_\eps'(\hat\varphi_\eps)\nabla\hat\varphi_\eps
\to m'(\hat\varphi)\nabla\hat\varphi \quad\text{a.e.~in } Q\,.
\]
Noting also that
\[
 |m_\eps'(\hat\varphi_\eps)\nabla\hat\varphi_\eps|^2\leq
 \norm{m'}_{L^\infty(-1,1)}|\nabla\hat\varphi_\eps|^2
\]
where the right-hand side is uniformly integrable $Q$
(because it converges strongly in $L^1(Q)$),
by Vitali's generalized dominated convergence theorem we infer that 
\[
  m_\eps'(\hat\varphi_\eps)\nabla\hat\varphi_\eps\to
  m'(\hat\varphi)\nabla\hat\varphi \quad\text{in } L^2(Q) \quad\hat\P\text{-a.s.}
\]
Since $\zeta\in L^\infty(Q)$ by choice of $\bar m$, this implies that 
\beq
  \label{conv_aux3}
  m'_\eps(\hat\varphi_\eps)\nabla\hat\varphi_\eps\cdot\eta\to
  m'(\hat\varphi)m'(\hat\varphi)\nabla\hat\varphi\cdot\eta
  \quad\text{in } L^2(0,T; H) \quad\hat\P\text{-a.s.}
\eeq
Let us focus now on the second term on the right-hand side of \eqref{def_eta_ep}.
Note that since $|\hat\varphi|\leq1$ 
and $\hat\varphi_\eps\to \hat\varphi$ almost everywhere in $Q$,
we have that (see \cite[pp.~416--417]{ell-gar})
\[
  m_\eps(\hat\varphi_\eps)F_\eps''(\hat\varphi_\eps)\to
  m(\hat\varphi)F''(\hat\varphi) \quad\text{a.e.~in } \hat\Omega\times Q\,.
\]
Since also $\nabla\hat\varphi_\eps\to\nabla\hat\varphi$
almost everywhere in $Q$,
recalling that $|m_\eps F_\eps''|\leq\norm{mF''}_{C^0([-1,1])}$ by 
{\bf D2}, the Vitali dominated convergence theorem yields again
\beq\label{conv_aux4}
  m_\eps(\hat\varphi_\eps)F_\eps''(\hat\varphi_\eps)\nabla\hat\varphi_\eps\to
  m(\hat\varphi)F''(\hat\varphi)\nabla\hat\varphi \quad\text{in } L^2(0,T; H)\,, \quad\hat\P\text{-a.s.}
\eeq
Consequently, using \eqref{conv_aux1}, \eqref{conv_aux2}, \eqref{conv_aux3}, and \eqref{conv_aux4}
in \eqref{def_eta_ep} and letting $\eps\to0$, we obtain exactly
the variational formulation \eqref{def_eta} for every $\zeta\in L^\infty(0,T; V_2\cap H^{2\bar m}(\OO))$.
A classical density argument yields then
\eqref{def_eta} for every $\zeta\in L^\infty(0,T; V_1\cap L^d(\OO))$.
Finally, the variational formulation \eqref{var-for-D} follows then
from \eqref{var_eta} and \eqref{def_eta} with the choice $\zeta=\nabla v$.
The regularity $\hat\varphi\in C^0_w([0,T]; L^2(\hat\Omega; V_1))$ is 
a consequence of the regularity
\[
\hat\varphi\in L^\infty(0,T; L^2(\hat\Omega; V_1))\cap 
C^0([0,T]; L^2(\Omega; (V_2\cap W^{1,d}(\OO))^*)\,,
\]
obtained by comparison in \eqref{var-for-D}.

The only thing that remains to be proved is the
regularity of $F(\hat\varphi)$ and $M(\hat\varphi)$.
Since we have the convergence
$\hat\varphi_\eps(t)\to\hat\varphi(t)$ almost everywhere
in $\hat\Omega\times \OO$, there is a measurable
set $A_0\subset \hat\Omega\times \OO$
(possibly depending on $t$) of full measure such that 
$\hat\varphi_\eps(t)\to\hat\varphi(t)$ pointwise in $A_0$ for every $t\in[0,T]$.
Let us show that 
\beq\label{liminf_eps}
  F_1(\hat\varphi(t))\leq\liminf_{\eps\to0}F_{1,\eps}(\hat\varphi_\eps(t))\,, 
  \quad
  M(\hat\varphi(t))\leq\liminf_{\eps\to0}M_\eps(\hat\varphi_\eps(t))
  \qquad\text{in } A_0 \quad\forall\,t\in[0,T]\,.
\eeq
Indeed, since we already know that 
$|\hat\varphi|\leq1$, for every
$(\hat\omega,x)\in A_0$ it holds either
$|\hat\varphi(\hat\omega,t,x)|<1$ or $|\hat\varphi(\hat\omega,t,x)|=1$.
In the former case, since $\hat\varphi_\eps(\hat\omega,t,x)\to\hat\varphi(\hat\omega,t,x)$,
we have $|\hat\varphi_\eps(\hat\omega,t,x)|<1$ for $\eps$ small enough, 
hence also $F_{1,\eps}(\hat\varphi_\eps(\hat\omega,t,x))=F(\hat\varphi_\eps(\hat\omega,t,x))$
and $M_{\eps}(\hat\varphi_\eps(\hat\omega,t,x))=M(\hat\varphi_\eps(\hat\omega,t,x))$
for $\eps$ small enough, from which \eqref{liminf_eps} follows
by continuity of $F_{1}$ and $M$ in $(-1,1)$.
In the latter case, it holds either $\hat\varphi(\hat\omega,t,x)=1$
or $\hat\varphi(\hat\omega,t,x)=-1$.
In the first possibility, for $\eps$ small we have
by definition of $F_{1,\eps}$ 
\[
  F_{1,\eps}(\hat\varphi_\eps(\hat\omega,t,x))\geq
  \min\{F_1(\hat\varphi_\eps(\hat\omega,t,x)), F_1(1-\eps)\}\to L\,,
\]
where $L:=\lim_{r\to1^-}F_1(r)\in[0,+\infty]$ is well defined by monotonicity of $F_1$.
If $L=+\infty$, then by comparison we deduce that 
$F_{1,\eps}(\hat\varphi_\eps(\hat\omega,t,x))\to+\infty$, so that 
\eqref{liminf_eps} holds automatically. If $L<+\infty$, then 
$F_1$ can be extended by continuity in $1$, and
the inequality gives
\[
  \liminf_{\eps\to0}F_{1,\eps}(\hat\varphi_\eps(\hat\omega,t,x))\geq
  L=F_1(1)=F_1(\hat\varphi(\hat\omega,t,x))\,,
\]
so that \eqref{liminf_eps} still holds.
The argument for the point $-1$
and the function $M$ is exactly the same, and 
\eqref{liminf_eps} is proved.
Finally, the inequalities \eqref{liminf_eps}, \eqref{est2_eps},
and \eqref{est6_eps}, together with the Fatou's lemma,
yield the desired assertion. This completes the proof of Theorem~\ref{thm2}.

\section*{Acknowledgments}
The author gratefully acknowledges financial support from the Austrian Science Fund (FWF)
through project M 2876. 
The author is also very grateful to the anonymous referees for their 
valuable comments and constructive suggestions.

%%%%%%%%%%%%%%%%%%%%%%%%%%%%%%%%%%%%%%%%

%\bibliography{ref}{}
\bibliographystyle{abbrv}

\def\cprime{$'$}

\end{document}